\documentclass[twoside,leqno]{article}

\usepackage[letterpaper]{geometry}

\usepackage{siamproceedings}

\usepackage[T1]{fontenc}
\usepackage{amsfonts}
\usepackage{graphicx}
\usepackage{epstopdf}
\usepackage{enumitem}
\usepackage{algorithmic}
\ifpdf
  \DeclareGraphicsExtensions{.eps,.pdf,.png,.jpg}
\else
  \DeclareGraphicsExtensions{.eps}
\fi

\usepackage{sfilip_abbreviations}

\crefname{subsection}{\S}{\S}

\newif\ifdebug
\debugfalse

\ifdebug

	\usepackage[layout = {nomargin, inline, index}]{fixme}
	\fxsetup{marginclue, theme = color, mode = multiuser, draft}
	\FXRegisterAuthor{sf}{esf}{\color{blue}Simion}

	\RequirePackage[orthodox]{nag}

\else
	\usepackage{fixme}
	\fxsetup{final}
\fi

\newsiamremark{remark}{Remark}

\newsiamremark{example}{Example}
\crefname{example}{Example}{Examples}

\newsiamremark{hypothesis}{Hypothesis}
\crefname{hypothesis}{Hypothesis}{Hypotheses}

\newsiamthm{conjecture}{Conjecture}
\crefname{conjecture}{Conjecture}{Conjectures}

\newsiamthm{problem}{Problem}
\crefname{problem}{Problem}{Problems}

\newsiamthm{claim}{Claim}

\usepackage{amsopn}

\begin{document}

\newcommand\relatedversion{}
\renewcommand\relatedversion{\thanks{The full version of the paper can be accessed at \protect\url{https://arxiv.org/abs/0000.00000}}} % Replace URL with link to full paper or comment out this line

\title{Measure Rigidity beyond Homogeneous Dynamics}
    \author{Simion Filip
    \thanks{
        University of Chicago (\email{sfilip@math.uchicago.edu}),
        Revised \today}
        }

\date{October 2025}

\maketitle

% Copyright Statement
% When submitting your final paper to a SIAM proceedings, it is requested that you include
% the appropriate copyright in the footer of the paper.  The copyright added should be
% consistent with the copyright selected on the copyright form submitted with the paper.
% Please note that "20XX" should be changed to the year of the meeting.

% Default Copyright Statement
\fancyfoot[R]{\scriptsize{Copyright \textcopyright\ 20XX by SIAM\\
Unauthorized reproduction of this article is prohibited}}

% Depending on which copyright you agree to when you sign the copyright form, the copyright
% can be changed to one of the following after commenting out the default copyright statement
% above.

%\fancyfoot[R]{\scriptsize{Copyright \textcopyright\ 20XX\\
%Copyright for this paper is retained by authors}}

%\fancyfoot[R]{\scriptsize{Copyright \textcopyright\ 20XX\\
%Copyright retained by principal author's organization}}

%\pagenumbering{arabic}
%\setcounter{page}{1}%Leave this line commented out.

\begin{abstract}
    We describe recent work that extends some of the measure and topological rigidity results in dynamical systems from situations homogeneous under a Lie group to quite general manifolds.
\end{abstract}

% \listoffixmes

% AMS MSC codes 37C40, 37D25, 32G15.

% keywords: Measure rigidity, random dynamical systems, generalized u-Gibbs states, normal forms, homogeneous dynamics, Anosov flows, character varieties, K3 surfaces, symplectic dynamics

%%%%%%%%%%%%%%%%%%%%%%%%%%%%%%%%%%%%%%%%%%%%%%%%%%%%%%%%%%%%%%%%%%%%%%%%%%%%%%%
%%%%%%%%%%%%%%%%%%%% Section Introduction %%%%%%%%%%%%%%%%%%%%
%%%%%%%%%%%%%%%%%%%%%%%%%%%%%%%%%%%%%%%%%%%%%%%%%%%%%%%%%%%%%%%%%%%%%%%%%%%%%%%

                    \section{Introduction.}
                    \label{sec:Introduction}

For most dynamical systems the classification of orbits and their statistical behavior, i.e. the invariant measures, is an intractable problem.
There exist however interesting dynamical systems for which such a classification is possible, or at least conjectured, and one speaks of ``rigidity'' of such systems.
A first instance is a theorem of Kronecker, which says that for an irrational translation on $\bR/\bZ$, \emph{every} orbit is dense, and a theorem of Weyl guarantees that the orbit is in fact equidistributed with respect to Lebesgue measure.

This has been generalized to dynamical systems on homogeneous spaces of Lie groups, as we will recall in more detail below.
This paper describes an extension to certain classes of dynamical systems on general manifolds, without assuming any homogeneity.

\noindent\textbf{Rigidity in Homogeneous Dynamics.}
The proofs of Kronecker's and Weyl's theorems extend almost verbatim to translations on a compact abelian group (or any compact group for that matter), but the orbit closures are now given by cosets of closed subgroups.
The next class of systems for which analogous results are known is that of unipotent flows on finite-volume quotients of Lie groups -- these are Raghunathan's conjectures established in full generality by Ratner \cite{Ratner1990_On-measure-rigidity-of-unipotent-subgroups-of-semisimple-groups,Ratner1991_On-Raghunathans-measure-conjecture,Ratner1991_Raghunathans-topological-conjecture-and-distributions-of-unipotent-flows} with important contributions of Margulis \cite{Margulis1987_Formes-quadratriques-indefinies-et-flots-unipotents-sur-les-espaces-homogenes}, Dani--Margulis \cite{MR1101994,DaniMargulis1993_Limit-distributions-of-orbits-of-unipotent-flows-and-values}, Margulis--Tomanov \cite{MargulisTomanov1994_Invariant-measures-for-actions-of-unipotent-groups-over-local}.
Extra invariance by unipotent elements is key to many rigidity results, and as we shall see below these arise in all systems with some hyperbolicity, even if they are not postulated at the start.
Another source of rigidity comes from higher rank features of Lie groups and manifests itself in Margulis super-rigidity \cite{Margulis_Discrete-subgroups-of-semisimple-Lie-groups1991}, the Zimmer program \cite{Zimmer1987_Actions-of-semisimple-groups-and-discrete-subgroups} (see Brown \cite{Brown2023_Lattice-subgroups-acting-on-manifolds} and Fisher \cite{Fisher2023_Rigidity-lattices-and-invariant-measures-beyond-homogeneous} for the most recent developments), and rigidity of abelian actions (see Lindestrauss's report \cite{Lindenstrauss2010_Equidistribution-in-homogeneous-spaces-and-number-theory}).
The underlying mechanisms behind this type of rigidity are however of a different flavor than those discussed in this text.

\noindent\textbf{Rigidity for Random Walks.}
A substantial advance was obtained in the work of Benoist--Quint \cite{BenoistQuint2011_Mesures-stationnaires-et-fermes-invariants-des-espaces-homogenes,BenoistQuint2013_Stationary-measures-and-invariant-subsets-of-homogeneous-spaces-II}, who obtained orbit closure and measure rigidity results analogous to Ratner's theorems, but for random walks on Lie groups.
Crucially, the dynamics can be generated by exclusively hyperbolic elements, but ultimately extra invariance is obtained along unipotent directions.
This was extended by Eskin--Mirzakhani--Mohammadi \cite{EskinMirzakhaniMohammadi2015_Isolation-equidistribution-and-orbit-closures-for-the-rm-SL2R-action-on-moduli,EskinMirzakhani2018_Invariant-and-stationary-measures-for-the-rm-SL2R-action-on-moduli-space} to the action of $\SL_2(\bR)$ on strata of translation surfaces, and by Brown--Rodriguez Hertz \cite{BrownRodriguez-Hertz2017_Measure-rigidity-for-random-dynamics-on-surfaces-and-related-skew} to random walk on real $2$-dimensional surfaces, see also Cantat--Dujardin \cite{CantatDujardin2023_Random-dynamics-on-real-and-complex-projective-surfaces} for the complex $2$-dimensional case, and Eskin--Lindenstrauss \cite{EskinLindenstrauss2018_Random-walks-on-locally-homogeneous-spaces} for homogeneous random walks.

\noindent\textbf{The general case.}
The theme of this text is that it is possible to extend some of the above rigidity results to a fairly general situation.
One class of examples comes from random walks generated by diffeomorphisms of smooth manifolds, and another from general $1$-parameter smooth maps (with discrete or continuous time) that have some initial regularity as codified by the notion of ``generalized u-Gibbs state'' (\cref{def:Generalized-u-Gibbs-state}).
There is a known analogy between the two situations.

The main technical result, \cref{thm:Inductive-step} below, shows that the relevant stationary or invariant measures must have extra invariance, unless their stable and unstable supports satisfy a natural ``joint integrability'' condition.
In situations where one has further control on the hyperbolicity, e.g. through the notion of uniform expansion, this information can be used to deduce further properties of the measures.

This text can serve as an introduction to the techniques and results from the author's joint work with Brown, Eskin, and Rodriguez Hertz \cite{BrownEskinFilip2025_Measure-rigidity-for-generalized-u-Gibbs-states-and-stationary}.
We include in the discussion a motivating rigidity result for totally geodesic manifolds on analytic Riemannian manifolds with variable negative sectional curvature obtained with Fisher and Lowe \cite{FilipFisherLowe2024_Finiteness-of-totally-geodesic-hypersurfaces}.

\noindent\textbf{Overview of Contents.}
We begin with some basic notation and three results that are elementary to state in \cref{sec:background_and_results}.
Next, we introduce the background from smooth dynamics, such as stable/unstable manifolds, and entropy theory, in \cref{sec:Tools-from-Smooth-Dynamics}.
The key notion of normal forms is introduced in \cref{sec:Normal-Forms}.
It is the tool that allows us to express unipotent invariance of conditional measures along stable or unstable manifolds.
The associated linearization cocycles are, among many other things, used to track divergence of the geometric objects along which one has extra invariance.
We then state the main technical result, and indicate some of the ingredients of the proof, in \cref{sec:Extra-Invariance-or-Joint-Integrability}.
We end in \cref{sec:Further-Directions} with a discussion of some natural questions that are motivated by the results we discuss.

%%%%%%%%%%%%%%%%%%%%%%%%%%%%%%%%%%%%%%%%%%%%%%%%%%%%%%%%%%%%%%%%%%%%%%%%%%%%%%%%
%%%%%%%%%%%% End of Section Introduction %%%%%%%%%%%%%%%%%%%%
%%%%%%%%%%%%%%%%%%%%%%%%%%%%%%%%%%%%%%%%%%%%%%%%%%%%%%%%%%%%%%%%%%%%%%%%%%%%%%%%

%%%%%%%%%%%%%%%%%%%%%%%%%%%%%%%%%%%%%%%%%%%%%%%%%%%%%%%%%%%%%%%%%%%%%%%%%%%%%%%
%%%                 Start of Section: Background and Results
%%%%%%%%%%%%%%%%%%%%%%%%%%%%%%%%%%%%%%%%%%%%%%%%%%%%%%%%%%%%%%%%%%%%%%%%%%%%%%%

\section{Background and Results.}
    \label{sec:background_and_results}
We introduce the context and some notation for the dynamical systems that we consider in \cref{ssec:Setup-Manifolds-Flows-Random-Walks}.
Three results that motivate the more technical main theorems are presented in \cref{ssec:Some-sample-theorems}.

%%%%%%%%%%%%%%%%%%%%%%%%%%%%%%%%%%%%%%%%%%%%%%%%%%%%%%%%%%%%%%%%%%%%%%%%%%%%%%%%
%%%%% SubSection Setup: Manifolds, Flows, Random Walks %%%%%

    \subsection{Setup: Manifolds, Flows, Random Walks.}
        \label{ssec:Setup-Manifolds-Flows-Random-Walks}
\hfill

\noindent \textbf{Standing Assumptions.}
Let $Q$ be a manifold and denote by $\Diff(Q)$ the group of diffeomorphisms of $Q$.
We can and will specialize later to diffeomorphisms preserving additional structures on $Q$, such as a smooth volume, a symplectic form, a real or complex analytic structure, etc.
We also note here that $Q$ is not required to be compact, and it might not even be complete with respect to a natural underlying Riemannian metric: for example, it could be the smooth locus of a singular algebraic variety such as the character variety of a surface group representation.

Our main result \cref{thm:Inductive-step} applies quite generally to a flow on a manifold with some hyperbolicity.
We will aim for notational simplicity below and state results in the setting where they are easiest to formulate.
Translations between random walk and flow/group action settings are possible but occasionally omitted in this text.

\noindent \textbf{Random Walks.}
Suppose that $\mu$ is a probability measure on $\Diff(Q)$; precise technical assumptions on $Q$ and $\mu$ are spelled out in \cite[\S2]{BrownEskinFilip2025_Measure-rigidity-for-generalized-u-Gibbs-states-and-stationary}.
A typical example arises from the action of a discrete group $\Gamma$ on $Q$:  define $\mu$ by assigning nonzero probabilities to the generators of $\Gamma$ and consider the resulting random walk on $Q$.

Recall next that a probability measure $\nu$ on $Q$ is called \emph{$\mu$-stationary} if $\nu = \mu * \nu$ where
\[
    \mu * \nu := \int_{\Diff(Q)} f_* \nu \, d\mu(f) \text{ and } f_* \nu(A) := \nu\left(f^{-1}(A)\right).
\]
The convolution $\mu *\mu$ is defined similarly and gives again a probability measure on $\Diff(Q)$, and we write $\mu^{*n}$ for the $n$-th convolution power.
We restrict to ergodic stationary measures, i.e. those that cannot be expressed as a nontrivial convex combination of other stationary measures.
Recall also that the Krylov--Bogolyubov argument shows stationary measures always exist (when $Q$ is compact).

A first manifestation of measure rigidity occurs when a measure is not just stationary but invariant under the elements in the support of $\mu$.
For example, the uniform measure on the middle third Cantor set is stationary under the action of the two maps $x\mapsto x/3$ and $x\mapsto (x+2)/3$ (with probability $1/2$ each), but it is not invariant under either of the two maps.

\noindent \textbf{Flows.}
We will also consider smooth $1$-parameter flows $\{g_t\}_{t\in \bR}$ acting on $Q$, with $t\mapsto g_t$ a group homomorphism from $(\bR,+)$ to $\Diff(Q)$.
In this case, $\nu$ is $g_t$-invariant if $(g_t)_* \nu = \nu$ for all $t\in \bR$ and ergodicity is defined analogously to the random walk case.

\noindent \textbf{Suspension and Skew-Products.}
A general construction allows us to change setting when needed.

\noindent \textit{Discrete to Continuous.}
First, given a single smooth diffeomorphism $f\in \Diff(Q)$, we can consider the suspension space $Q' := (Q\times \bR)/\bZ$ where the $\bZ$-action is generated by $(x,t)\mapsto (f(x),t+1)$ and define the flow $g_t(x,s) := (x,s+t)$.
This reduces discrete time examples to continuous time ones.

More generally, if $\Gamma$ is a lattice in a Lie group $G$, and $\Gamma$ acts on the manifold $Q$, we can consider the suspension space $Q' := \leftquot{\Gamma}{G\times Q}$ where $\Gamma$ acts diagonally by $\gamma\cdot (g,x) := (\gamma g,\gamma x)$.
We then have a $G$-action on $Q'$ by right multiplication on the first factor, and we can take $1$-parameter subgroups $g_t\subset G$ to obtain flows on $Q'$.
A useful and nontrivial example of the above setup is for $\Gamma$ to be a free group, which embeds as a lattice in $G=\SL_2(\bR)$.

\noindent \textit{Bernoulli shifts.}
Finally, we can consider suspensions over Bernoulli shifts.
Set $B^+:= (\Diff(Q))^{\bZ_{\geq 0}}$ and let $F\colon B^+\to B^+$ be the left shift map; the measure $\widehat{\mu}^+:=\mu^{\otimes \bZ_{\geq 0}}$ is $F$-invariant and ergodic.
We will write $b=(b_0,b_1,\ldots)$ for elements of $B^+$.
We can next consider the space $\widehat{Q}^+:=Q\times B^+$, with skew-product map $\hat{F}(q,b):=(b_0(q), F(b) )$ and measure $\widehat{\nu}^+ := \widehat{\mu}^+\otimes \nu$, which is $F$-invariant if and only if $\nu$ is $\mu$-stationary.

It is useful to also consider $B:=\Diff(Q)^{\bZ}$ and associated $\widehat{Q}$, with shift and skew product maps defined in the same way.
The martingale convergence theorem constructs a natural measure $\widehat{\nu}$ associated to the $\mu$-stationary measure $\nu$.
For $b\in B$ we will write $b_{\geq 0}\in B^+$ for its components with index in $\bZ_{\geq 0}$.
Additionally, for $b=(\dots,b_{-1}, b_0, b_1,\dots)$ we will write $b^{\circ n}:=b_{n}\circ\dots\circ b_0$ for the composed diffeomorphism.
See \cite[\S2]{BrownEskinFilip2025_Measure-rigidity-for-generalized-u-Gibbs-states-and-stationary} for more on these constructions.

% \fxnote{Include a dictionary between the two situations, Eskin--Margulis nondivergence, u-Gibbs states and stationary measures. Is it really three situations? Unipotent Flows (or P-invariant measures)?, Random Walks, u-Gibbs states}

%%%%% End of SubSection Setup: Manifolds, Flows, Random Walks %%%%%
%%%%%%%%%%%%%%%%%%%%%%%%%%%%%%%%%%%%%%%%%%%%%%%%%%%%%%%%%%%%%%%%%%%%%%%%%%%%%%%%

%%%%%%%%%%%%%%%%%%%%%%%%%%%%%%%%%%%%%%%%%%%%%%%%%%%%%%%%%%%%%%%%%%%%%%%%%%%%%%%%
%%%%% SubSection Three sample theorems %%%%%

    \subsection{Some sample theorems.}
        \label{ssec:Some-sample-theorems}
We start with several elementary statements that nonetheless require our main technical results.
A convenient assumption in the statements, which is made explicit below in \cref{def:Uniform-Expansion-and-Uniform-Gaps}, is that of ``uniform expansion'' and ``uniform gaps'':

\begin{theorem}[Finite or Uniform]
    \label{thm:Finite-or-Uniform}
    Suppose that $Q$ is compact and that the support of $\mu$ in $\Diff(Q)$ is finite and every element preserves a fixed smooth volume form.
    Suppose furthermore that $\mu$ is uniformly expanding and has uniform gaps for every $d\in \{1, \ldots , \dim(Q)-1\}$.
    
    Then every ergodic $\mu$-stationary measure on $Q$ is either finitely supported or is the smooth volume.
\end{theorem}
The case when $Q=\bR^d/\bZ^d$ and $\mu$ is supported on elements of $\SL_d(\bZ)$ was treated by Bourgain--Furman--Lindenstrauss--Mozes \cite{BourgainFurmanLindenstrauss2011_Stationary-measures-and-equidistribution-for-orbits-of-nonabelian-semigroups-on-the-torus} by Fourier-theoretic techniques, and Benoist--Quint \cite{BenoistQuint2011_Mesures-stationnaires-et-fermes-invariants-des-espaces-homogenes} with methods closer to our work.

\begin{theorem}[Stiffness]
    \label{thm:stiffness}
    Suppose that $Q$ is compact and that the support of $\mu$ on $\Diff(Q)$ is finite.
    Let $\nu$ be a $\mu$-stationary measure with the property that for every $\nu$-a.e. defined subbundle $S\subseteq TQ$ which is invariant by the elements in the support of $\mu$, the sum of Lyapunov exponents is nonnegative.

    Then $\nu$ is $\mu$-invariant.
\end{theorem}
The assumption above on finite support of $\mu$ is technical and should be easily removed once additional tools in the entropy theory of skew products is available.
Technically, the entropy of a continuous measure on $\Diff(Q)$ is infinite and this obstructs some entropy calculations as in \cref{ssec:Scheme-of-Proof}, but this is not essential.

The assumption in \cref{thm:stiffness} on Lyapunov exponents, for every $\mu$-stationary measure, holds if we have uniform expansion as in the next definition:

\begin{definition}[Uniform Expansion and Uniform Gaps]
    \label{def:Uniform-Expansion-and-Uniform-Gaps}
    For $d\in \{1, \ldots , \dim(Q)-1\}$, define $\Gr_d(TQ)$ to be the Grassmannian bundle of $d$-planes in the tangent bundle $TQ$.
    Identify a $d$-plane $P\in \Gr_d(T_q Q)$ with a vector $\xi_{P}\in \Lambda^d(T_qQ)$ and define
    \[
        \sigma(\xi,N):=\int_{\Diff(Q)} \log \left(\frac{\norm{D_q f^N(\xi)}}{\norm{\xi}}\right) \, d\mu^{*N}(f).
    \]
    The measure $\mu$ is \emph{uniformly expanding in dimension $d$} if there exists $N\geq 1$ and $c>0$ such that for every $P\in \Gr_d(TQ)$ we have
    \[
        \sigma(\xi_P,N) \geq c.
    \]
    Similarly, $\mu$ has \emph{uniform gaps in dimension $d$} if there exists $N\geq 1$ and $c>0$, as well as $\delta_d\in\{\pm 1\}$, such that for every $P_i\in \Gr_{d+i}(TQ)$, for $i\in \{0,1\}$ and with $P_0\subset P_1$, we have that
    \[
        \delta_d \cdot \left(\sigma(\xi_{P_1},N) - \sigma(\xi_{P_0},N)\right) \geq c.
    \]
\end{definition}
Informally, uniform expansion says that \emph{every} $d$-plane is expanded on average, while uniform gaps says that the $d$-th Lyapunov exponent (see \cref{ssec:Lyapunov-exponents}) is uniformly separated from zero.
It can be checked that these properties are independent of the choice of Riemannian metric, as long as the two Riemannian metrics are comparable up to a uniform constant.
We can extend the definition to refer to only a class of subspaces invariant under the dynamics, such as the isotropic ones in symplectic dynamics or complex ones in holomorphic dynamics.

We now discuss another result that extends rigidity results from the homogeneous to the inhomogeneous setting, obtained jointly with D.~Fisher and B.~Lowe \cite{FilipFisherLowe2024_Finiteness-of-totally-geodesic-hypersurfaces}:
\begin{theorem}[Totally Geodesic Rigidity]
    \label{thm:Totally-Geodesic-Rigidity}
    Suppose $M$ is a compact real-analytic Riemannian manifold with strictly negative sectional curvatures and $\dim(M)\geq 3$.
    If $M$ admits infinitely many distinct immersed closed totally geodesic hypersurfaces, then $M$ has constant negative curvature.
\end{theorem}
In the homogeneous setting, Bader--Fisher--Miller--Stover \cite{BaderFisherMiller2021_Arithmeticity-superrigidity-and-totally-geodesic-submanifolds} in general and Margulis--Mohammadi \cite{MohammadiMargulis2022_Arithmeticity-of-hyperbolic-3-manifolds-containing-infinitely-many-totally} for $n=3$ showed that a finite volume constant negative curvature manifold with infinitely many totally geodesic hypersurfaces must in fact be arithmetic.
We discuss some expectations in this direction in \cref{ssec:rigidity_and_special_orbit_closures} below.

We also note that the main result in \cite{BrownEskinFilip2025_Measure-rigidity-for-generalized-u-Gibbs-states-and-stationary} is used there to derive a generalization of \cite{EskinMirzakhani2018_Invariant-and-stationary-measures-for-the-rm-SL2R-action-on-moduli-space} to the case of the action of $\SL_2(\bR)$ on products of strata of translation surfaces.

%%%%% End of SubSection Some sample theorems %%%%%
%%%%%%%%%%%%%%%%%%%%%%%%%%%%%%%%%%%%%%%%%%%%%%%%%%%%%%%%%%%%%%%%%%%%%%%%%%%%%%%

%%%%%%%%%%%%%%%%%%%%%%%%%%%%%%%%%%%%%%%%%%%%%%%%%%%%%%%%%%%%%%%%%%%%%%%%%%%%%%%
%%%                 End of Section: Background and Results
%%%%%%%%%%%%%%%%%%%%%%%%%%%%%%%%%%%%%%%%%%%%%%%%%%%%%%%%%%%%%%%%%%%%%%%%%%%%%%%

%%%%%%%%%%%%%%%%%%%%%%%%%%%%%%%%%%%%%%%%%%%%%%%%%%%%%%%%%%%%%%%%%%%%%%%%%%%%%%%
%%%%%%%%%%%%%%%%%%%% Section Tools from Smooth Dynamics %%%%%%%%%%%%%%%%%%%%
%%%%%%%%%%%%%%%%%%%%%%%%%%%%%%%%%%%%%%%%%%%%%%%%%%%%%%%%%%%%%%%%%%%%%%%%%%%%%%%

                    \section{Tools from Smooth Dynamics.}
                    \label{sec:Tools-from-Smooth-Dynamics}

To provide more context and some ideas from the proofs, we recall first some basic results from the theory of smooth dynamical systems.

%%%%%%%%%%%%%%%%%%%%%%%%%%%%%%%%%%%%%%%%%%%%%%%%%%%%%%%%%%%%%%%%%%%%%%%%%%%%%%%%
%%%%% SubSection Lyapunov exponents %%%%%

    \subsection{Lyapunov exponents.}
        \label{ssec:Lyapunov-exponents}

\hfill

\noindent \textbf{Cocycles.}
Suppose, quite generally, that we have a space $Q$ with an action of a group or semigroup $\Gamma$.
A \emph{cocycle} over the $\Gamma$-action is the data of a vector bundle $E\to Q$ and a lift of the $\Gamma$-action to $E$ by vector bundle automorphisms.
Concretely, for each $\gamma\in \Gamma$ we have for each $q\in Q$ a linear map $\gamma_E\colon E(q)\to E(\gamma q)$, of the same regularity in $q$ as the vector bundle $E$, and obeying the natural compatibility conditions for the product of two group elements.
See \cref{ssec:A-Hierarchy-of-Cocycles} below for basic examples of smooth and measurable cocycles.

The fundamental result concerning the asymptotic behavior of cocycles is the Oseledets Multiplicative Ergodic Theorem (MET):

\begin{theorem}[MET]
    \label{thm:MET}
    Let $f\colon Q\to Q$ be a invertible measurable transformation, preserving an ergodic probability measure $\nu$.
    Let $E\to Q$ be a cocycle over $f$, equipped with a measurable norm $\|\cdot\|$ and satisfying the integrability condition $\int_Q \log^+ \left(\norm{f_E},\norm{f^{-1}_E}\right) \, d\nu < \infty$.

    Then there exist numbers $\lambda_1 > \lambda_2 > \cdots > \lambda_k$ (the Lyapunov exponents) and a measurable $f_E$-invariant splitting $E = E_1 \oplus \cdots \oplus E_k$ such that for $\nu$-almost every $q\in Q$ and every nonzero $v\in E_i(q)$ we have that
    \[
        \lim_{n\to \infty} \frac{1}{n} \log \|f_E^{\circ n}(v)\| = \lambda_i.
    \]
\end{theorem}
Let us emphasize that even if $E$ is a smooth vector bundle, the Oseledets splitting is only measurable in general.

\noindent \textbf{Jet bundles.}
Recall that $Q$ is a smooth manifold, and has in particular a tangent bundle $TQ$ which becomes a cocycle for any action of a group by diffeomorphisms.
More generally, it is important to consider the vector bundle $J^rQ$ of $r$-jets of functions on $Q$, and the related princial bundles $F^rQ$ of $r$-jets of charts on $Q$: the fiber of $F^rQ$ at $q\in Q$ consists of equivalence classes of local diffeomorphisms $\phi\colon (\bR^{\dim(Q)},0)\to (Q,q)$, where two such maps are equivalent if their Taylor expansions agree up to order $r$.

Although the Lyapunov spectrum of $J^rQ$ is determined by that of $TQ$, the cocycle $J^rQ$ is important because it encodes higher order information about the dynamics.
This plays a particularly important role in the theory of normal forms, see \cref{sec:Normal-Forms} below.

%%%%% End of SubSection Lyapunov exponents %%%%%
%%%%%%%%%%%%%%%%%%%%%%%%%%%%%%%%%%%%%%%%%%%%%%%%%%%%%%%%%%%%%%%%%%%%%%%%%%%%%%%%

%%%%%%%%%%%%%%%%%%%%%%%%%%%%%%%%%%%%%%%%%%%%%%%%%%%%%%%%%%%%%%%%%%%%%%%%%%%%%%%%
%%%%% SubSection Invariant Manifolds %%%%%

    \subsection{Invariant Manifolds.}
        \label{ssec:Invariant-Manifolds}

We continue working with a manifold $Q$ and a diffeomorphism $f\in \Diff(Q)$ (or flow $g_t$) preserving an ergodic probability measure $\nu$, or in the case of random walks a sequence of diffeomorphisms $f_n\in \Diff(Q)$ chosen independently at random with law $\mu$ and a $\mu$-stationary measure $\nu$.
We denote by $f^{\circ n}$ the $n$-th iterate of $f$, or $f^{\circ n} := f_{n-1}\circ \cdots \circ f_0$ in the case of random walks.
The Lyapunov exponents of the tangent bundle $TQ$ are crucial invariants of the dynamics, and we will assume that at least some are nonzero.
We then have a (measurable) decomposition of the tangent bundle into stable, unstable, and center subbundles
\[  TQ = W^s \oplus W^0 \oplus W^u \]
where $W^s$ (resp. $W^u$) is the sum of the Oseledets subspaces corresponding to negative (resp. positive) Lyapunov exponents, and $W^0$ corresponds to zero exponents.
For simplicty, we describe the structures below in the stable directions, but analogous statements hold, under time reversal, for the unstable directions.

\noindent \textbf{Stable Manifolds.}
Under suitable technical assumptions, the stable subbundle $W^s$ is integrable to a measurable family of (local) stable manifolds for $\nu$-almost every $q\in Q$:
\[
    \cW^s[q]=\left\lbrace x\in Q \colon  d_Q(f^{\circ n} (x),f^{\circ n} (q))\le C(x)e^{-\ve n}\right \rbrace
\]
for any sufficiently small $\ve>0$ and some $C(x)>0$.
This family of manifolds is equivariant for the dynamics.

In the random walk case, on the space $\hat{Q}=Q\times B$, we have both the ``combinatorial'' local stable
\[
    \frakW^s_{1}[q,b]:=\{(q,b')\in \hat{Q}\colon b'_{\geq 0 }=b_{\geq 0}\}
\]
and the ``smooth'' stable:
\[
    \cW^{s}[q,b]:=\{(q',b)\in \hat{Q} \colon d_{Q}(b^{\circ n}(q),b^{\circ n}(q')\leq C(q',b)e^{-\ve n})\}.
\]
We note that $\frakW^s_1$ is the pullback of an analogous measurable partition on $B$, which can be used to compute the entropy of $\left(B,\mu^{\otimes \bZ}\right)$ under the shift.
When $\mu$ is atomic, this equals $H(\mu)$.
For future purposes, we note here that
\begin{align}
    \label{eqn:relative_entropy_bound}
    h_{\hat{\nu}}\left(\hat{F}^{-1};\frakW^s_1\right) \leq H(\mu)
\end{align}
using Jensen's inequality and the definition of entropy relative to a measurable partition; equality holds if and only if $\nu$ is $\mu$-invariant.

\noindent \textbf{Center-stable manifolds.}
The situation with the bundle $W^{cs}:=W^s\oplus W^0$ is more subtle, and in general there need not be a family of manifolds tangent to $W^{cs}$, which is also equivariant for the dynamics for all times.
What does exist is a formal germ of a manifold tangent to $W^{cs}$, which is equivariant for the dynamics.
One can either choose a smooth representative of this germ, or fix an arbitrarily large regularity class $C^r$ and construct a family of $C^r$-manifolds that are invariant under the dynamics but for finite time only.
Both notions are useful for the proof of the main result.

%%%%% End of SubSection Invariant Manifolds %%%%%
%%%%%%%%%%%%%%%%%%%%%%%%%%%%%%%%%%%%%%%%%%%%%%%%%%%%%%%%%%%%%%%%%%%%%%%%%%%%%%%%

%%%%%%%%%%%%%%%%%%%%%%%%%%%%%%%%%%%%%%%%%%%%%%%%%%%%%%%%%%%%%%%%%%%%%%%%%%%%%%%%
%%%%% SubSection Conditional Measures %%%%%

    \subsection{Conditional Measures.}
        \label{ssec:Conditional-Measures}

It is possible to disintegrate the measure $\nu$ along the unstable manifolds, obtaining a family of conditional measures $\nu^u[q]$ supported on $\cW^u[q]$ for $\nu$-almost every $q\in Q$.
To do so requires introducing a measurable partition of $Q$ subordinate to the stable manifolds, and Rokhlin's disintegration theorem.
Let us introduce the notation and recall the construction, since it will be useful later when discussing generalized u-Gibbs states in \cref{ssec:Generalized-u-Gibbs-states}, see also \cite[A.5]{BrownEskinFilip2025_Measure-rigidity-for-generalized-u-Gibbs-states-and-stationary} for details.

There exists a measurable partition of $(Q,\nu)$, denoted $\frakB_0$, which is subordinated to the unstable foliation, i.e. for $\nu$-a.e. $q\in Q$ the atom $\frakB_0[q]$ is contained in $\cW^u[q]$ and contains an open neighborhood of $q$.
Furthermore $\frakB_0$ is increasing under the dynamics, i.e. $f^{-1}(\frakB_0[f(q)])\subseteq \frakB_0[q]$.
Recall also that the partition $\frakB_0$ is ``$\nu$-measurable'' if there exists a measurable map $\pi_{\frakB_0}$ to a standard Borel space $S$, such that $\frakB_0$ is the preimage of the Borel $\sigma$-algebra mod sets of $\nu$-measure zero.
The conditional measures $\nu^u_{\frakB_0}[q]$ are then the disintegration of $\nu$ along the fibers of $\pi_{\frakB_0}$.
By using the increasing property of $\frakB_0$, one glues these conditional measures to $\nu^u[q]$ on $\cW^u[q]$, which are equivariant under the dynamics up to a constant scaling factor, and normalized by $\nu^u[q](\frakB_0[q])=1$.

We denote by $\frakB_0^-$ a partition of the stables with the same properties.

\noindent\textbf{Entropy and Dimension.}
Invariance of a stationary measure $\nu$ can be obtained using an equality between past and future entropy, based on the Ledrappier--Young formula \cite{LedrappierYoung1985_The-metric-entropy-of-diffeomorphisms.-II.-Relations-between-entropy-exponents} which relates the entropy of the measure $\nu$ to the dimensions of the conditional measures $\nu^u[q]$ along unstable manifolds and the Lyapunov spectrum.
We recall it in the following form.
Suppose $E_1\subseteq E_2\subseteq W^u$ are measurable dynamics-equivariant subbundles of the unstable tangent bundle, with the property that (appropriately interpreted) the conditional measures of $\nu^u$ along $E_1$ are Lebesgue class, and the support of $\nu^u$ is contained in $E_2$.
Denote by $E_i^{\lambda_j}$ the corresponding Lyapunov subspace.
The the entropy $h_{\nu}(f)$ of $f$ with respect to $\nu$ satisfies:
\begin{align}
    \label{eqn:Ledrappier_Young_Variant}
    \sum_{\lambda_i>0}\lambda_i \cdot\dim E_1^{\lambda_i}
    \leq h_{\nu}(f)
    \leq
    \sum_{\lambda_i>0} \lambda_i\cdot  \dim E_2^{\lambda_i}
\end{align}
Equality holds in the second inequality if and only if the conditionals of $\nu^u$ along $E_2$ are Lebesgue class.

%%%%% End of SubSection Conditional Measures %%%%%
%%%%%%%%%%%%%%%%%%%%%%%%%%%%%%%%%%%%%%%%%%%%%%%%%%%%%%%%%%%%%%%%%%%%%%%%%%%%%%%%

%%%%%%%%%%%%%%%%%%%%%%%%%%%%%%%%%%%%%%%%%%%%%%%%%%%%%%%%%%%%%%%%%%%%%%%%%%%%%%%%
%%%%% SubSection A Hierarchy of Cocycles %%%%%

    \subsection{A Hierarchy of Cocycles.}
        \label{ssec:A-Hierarchy-of-Cocycles}

We end these preliminaries with a description of several constructions of cocycles and their a priori properties.

\noindent \textbf{Smooth and Natural Cocycles.}
Starting from the manifold $Q$ and smooth action on it, we have the jet bundles described in \cref{ssec:Lyapunov-exponents}, which have a natural smooth structure.
We can further consider the associated natural maps, as well as linear-algebraic constructions such as quotients, tensor products, exterior and symmetric powers, duals, etc.
We will refer to any cocycle obtained by this procedure as ``smooth and natural''.

\noindent \textbf{Forward and Backward Dynamically Defined (fdd, bdd) Cocycles.}
For a smooth and natural cocycle $E$, the forward Oseledets filtrations $E^{\leq \lambda}$ are smooth along stable manifolds.
Cocycles obtained by standard linear-algebraic operations starting from the Oseledets filtrations of smooth and natural cocycles will be called \emph{forward dynamically defined} (fdd) cocycles.
By reversing time, and using the reverse filtrations $E^{\geq \lambda}$ we obtain \emph{backward dynamically defined} (bdd) cocycles.

\noindent\textit{Note:} A given measurable isomorphism class of cocycle can have both an fdd and bdd structure, so part of the structure of an fdd cocycle is the way the cocycle is constructed.
For example, an Oseledets subspace $E^{\lambda}$ has an fdd structure when written as $E^{\leq \lambda}/E^{<\lambda}$, but a bdd structure when written as $E^{\geq \lambda}/E^{>\lambda}$.

\noindent\textbf{Admissible Cocycles.}
Finally, we will consider arbitrary cocycles which arise as dynamics-equivariant subcocycles of smooth and natural cocycles, as well as those derived using the standard linear-algebraic operations from them.
We refer to these as \emph{admissible} cocycles, and all cocycles we work with are at least admissible.
A class of examples to keep in mind are the cocycles arising when a single Lyapunov subspace (which is admissible) is put into Jordan normal form, using further subcocycles, see \cite[B.3]{BrownEskinFilip2025_Measure-rigidity-for-generalized-u-Gibbs-states-and-stationary}.

An important feature of any admissible cocycle $E$ is that it admits a ``measurable connection'' $P^{\pm}$ along stable and unstable manifolds.
Specifically, for a set of full measure $Q_0$, if $x,y\in Q_0$ and $y\in \cW^u[x]$ we have a map $P^{+}(x,y)\colon E(x)\to E(y)$, which is equivariant for the dynamics and intertwines the measurable algebraic hulls.
Similar maps exist along stables, and we will sometimes denote them by $P^{u/s}$ instead of $P^{+/-}$.

%%%%% End of SubSection A Hierarchy of Cocycles %%%%%
%%%%%%%%%%%%%%%%%%%%%%%%%%%%%%%%%%%%%%%%%%%%%%%%%%%%%%%%%%%%%%%%%%%%%%%%%%%%%%%%

%%%%%%%%%%%%%%%%%%%%%%%%%%%%%%%%%%%%%%%%%%%%%%%%%%%%%%%%%%%%%%%%%%%%%%%%%%%%%%%%
%%%%%%%%%%%% End of Section Tools from Smooth Dynamics %%%%%%%%%%%%%%%%%%%%
%%%%%%%%%%%%%%%%%%%%%%%%%%%%%%%%%%%%%%%%%%%%%%%%%%%%%%%%%%%%%%%%%%%%%%%%%%%%%%%%

%%%%%%%%%%%%%%%%%%%%%%%%%%%%%%%%%%%%%%%%%%%%%%%%%%%%%%%%%%%%%%%%%%%%%%%%%%%%%%%%
%%%%%%%%%%%%%%%%%%%% Section Normal Forms %%%%%%%%%%%%%%%%%%%%
%%%%%%%%%%%%%%%%%%%%%%%%%%%%%%%%%%%%%%%%%%%%%%%%%%%%%%%%%%%%%%%%%%%%%%%%%%%%%%%%

                    \section{Normal Forms.}
                    \label{sec:Normal-Forms}

In this section, we recall in \cref{ssec:Classical-Normal-Forms} some features of the theory of normal forms, with some enhancements that are useful for our arguments.
In \cref{ssec:Normal-Forms-as-G-X-structures} we present an alternative approach to the construction of normal forms, which we hope simplifies some aspects and is also technically useful in deriving additional estimates for normal form coordinates.
Finally, in \cref{ssec:Generalized-u-Gibbs-states} we introduce the notion of ``generalized u-Gibbs state'' which plays a key role in our statements and proofs.

%%%%%%%%%%%%%%%%%%%%%%%%%%%%%%%%%%%%%%%%%%%%%%%%%%%%%%%%%%%%%%%%%%%%%%%%%%%%%%%%
%%%%% SubSection Classical Normal Forms %%%%%

    \subsection{Classical Normal Forms.}
        \label{ssec:Classical-Normal-Forms}

The theory of normal forms goes back to the work of Poincar\'e and Sternberg \cite{Sternberg1957_Local-contractions-and-a-theorem-of-Poincare} for a diffeomorphism near a fixed point.
The literature on the topic is vast, and below we will refer to the work of Kalinin and Sadovskaya \cite{KalininSadovskaya2017_Normal-forms-for-non-uniform-contractions} as a convenient reference (see also Melnick \cite{Melnick_Non-stationary-smooth-geometric-structures-for-contracting-measurable}), which generalized to the nonuniformly contracting setting earlier works, especially those of Guysinsky and Katok \cite{GuysinskyKatok1998_Normal-forms-and-invariant-geometric-structures-for-dynamical-systems}.

We deduce some additional properties of the formalism of normal forms.
Particularly useful for our arguments are various linearization constructions associated to normal forms, allowing us to further simplify the subresonant polynomial dynamics to a linear one.

\noindent \textbf{Gradings and Subresonance Polynomials.}
We need some preliminaries regarding polynomial maps of filtered vector spaces.
Let then $V$ be a finite-dimensional real vector space, equipped with a filtration by subspaces $0\subsetneq V^{\leq -\lambda_1} \subsetneq V^{\leq -\lambda_2} \subsetneq \dots \subsetneq V^{\leq -\lambda_n} = V$ where $\lambda_1 > \lambda_2 > \cdots > \lambda_n>0$ are positive real numbers.
Various linear-algebraic constructions inherit a natural filtration from the one on $V$.
Quite generally, a linear map $A\colon V\to W$ between filtered vector spaces satisfies $A\in \Hom(V,W)^{\leq \delta}$ if and only if $A\left(V^{\leq \lambda}\right) \subset W^{\leq \lambda+\delta}$ for all $\lambda$.
For example, the space of polynomial functions is the symmetric algebra $\Sym^\bullet(V^*)$, and the filtration is defined by declaring that $\xi\in V^{*,\leq \lambda_i}$ if $\xi\left(V^{\leq -\lambda_i}\right)=0$ and taking symmetric powers.

We can similarly define a filtration on the space of polynomial maps $\Poly(V,W)$, by declaring that $F\in \Poly(V,W)^{\leq \delta}$ if and only if $F^*\xi \in \Sym^\bullet(V^*)^{\leq \lambda+\delta}$ for every $\xi\in W^{*, \leq \lambda}$ and every $\lambda$.
Here $F^*\xi$ is the pullback of the linear function $\xi$ by the polynomial map $F$.

It is not hard to see that if all the weights $-\lambda_i$ that occur in the filtration of $V$ are strictly negative, then the space $\Poly(V,V)^{\leq 0}$ is finite-dimensional and closed under composition.

\begin{definition}[Subresonant Maps]
    \label{def:Subresonant-Maps}
    Define $\bbG^{sr}(V)$ to be the group of polynomial maps $F\in \Poly(V,V)^{\leq 0}$ that admit a polynomial inverse $F^{-1}\in \Poly(V,V)^{\leq 0}$ and call such maps \emph{subresonant}.

    Define $\bbG^{ssr}(V)$ to be the subgroup of $\bbG^{sr}(V)$ consisting of maps $F$ such that $D_0 F = \id$ and call such maps \emph{strictly subresonant}; equivalently, these are the maps $F$ such that $F - \id_V \in \Poly(V,V)^{< 0}$.
\end{definition}
For example, translations of $V$ are strictly subresonant maps.

The main result of \cite{KalininSadovskaya2017_Normal-forms-for-non-uniform-contractions} is that by a change of parametrization, the dynamics along stable manifolds can be put in normal form, i.e. to belong to the group $\bbG^{sr}(V)$.

\begin{theorem}[Normal Forms along Stables]
    \label{thm:Normal-Forms-along-Stables}
    Let $f\colon Q\to Q$ be a diffeomorphism preserving an ergodic probability measure $\nu$, and with additional integrability conditions on the higher jet cocycles.
    There exists a set $Q_0$ of full $\nu$-measure, such that for every $q\in Q_0$ there exists a $C^\infty$-diffeomorphism $H_q\colon \cW^s[q]\to W^s(q)$, where $W^s(q)\subset T_qQ$ is the subspace with negative Lyapunov exponents equipped with its Oseledets filtration, such that:
    \begin{enumerate}
        \item The maps $H_{f(q)}\circ f \circ H_q^{-1}\colon W^s(q)\to W^s(f(q))$ are subresonant polynomial maps.
        \item If $q'\in Q_0$ belongs to the local stable manifold $\cW^s[q]$, then the map $H_{q'}\circ H_q^{-1}\colon W^s(q)\to W^s(q')$ is a subresonant polynomial map.
    \end{enumerate}
\end{theorem}
The maps $H_q$ are called \emph{normal form coordinates}, depend measurably on $q$, but are not unique: if $G_q\in \bbG^{sr}(W^s(q))$ is a measurable family of subresonant polynomial maps, then $\widetilde{H}_q := G_q \circ H_q$ are also normal form coordinates.

\begin{remark}[Some further desiderata for normal forms]
    \label{rmk:Some-further-desiderata-for-normal-forms}
    First, it would be desirable to pin down a more intrinsic parametrization of stable manifolds, not up to the finite-dimensional ambiguity of precomposing by subresonant polynomial maps.
    For this, a space different from $W^s(q)$ is needed, one that includes higher order information about the dynamics.

    Second, it is desirable to have a \Holder-continuous dependence of the normal form coordinates on the basepoint $q$ (on sets of large measure).
    This is manifestly related to uniqueness.

    Both of these issues are addressed in \cref{ssec:Normal-Forms-as-G-X-structures} below, where a different approach to the construction of normal forms is presented.
    Moreover, as it will become apparent, this shows that the construction of normal forms, assuming the existence of stable manifolds, is in fact quite soft and involves integrating a ``smooth'' distribution, unlike the construction of the stable manifolds themselves.
\end{remark}

\noindent\textbf{Measurable connection for normal forms.}
The construction of a measurable connection on an admissible cocycle, discussed in \cref{ssec:A-Hierarchy-of-Cocycles}, applied to an appropriate cocycle constructed out of jets, gives a measurable family of maps $\wp^u(x,y)\in \bbG^{ssr}(\cW^u[x])$, for $x,y$ in a set of full measure.
These have the property that $\wp^u(x,y)(x)=y$ and are essential in establishing extra invariance of conditional measures.

\noindent\textbf{Normal forms for cocycles.}
An important technical ingredient is the construction of normal forms for cocycles, for example along unstables for cocycles smooth along unstables.
Just like normal forms for unstable manifolds bring the dynamics to a particularly simple polynomial form, cocycle normal forms bring any cocycle to a simpler polynomial form as well.

\noindent\textbf{An example of normal forms.}
Consider first $\bR^2$ with coordinate functions $(x,y)$ and respective weights $-\lambda_1<-\lambda_2<0$.
The most general subresonant maps are then of the following form, with linearization given by taking the Veronese embedding (with $k_0$ maximal such that $k_0\lambda_2\leq \lambda_1$ and assuming for simplicity $a=b=0$):

\begin{align*}
    \everymath{\displaystyle}
    \begin{bmatrix}
        x\\
        y
    \end{bmatrix}
    \xrightarrow{f}
    \begin{bmatrix}
        a + \lambda x + \sum_{0<k\lambda_2 \leq \lambda_1}c_k y^k\\
        b + \mu y \phantom{+ \sum_{0<k\lambda_2 \leq \lambda_1}c_k y^k}
    \end{bmatrix}
    \text{ and linearization }
     \underbrace{
    \begin{bmatrix}
        x\\
        y
    \end{bmatrix}
    }_{q}
    \leftrightarrow
    \underbrace{
    \begin{bmatrix}
        y^{k_0}\\
        \vdots\\
        y\\
        x
    \end{bmatrix}
    }_{Lq}
    \mapsto
    \underbrace{
    \begin{bmatrix}
        \mu^{k_0} & 0 & \hdots & 0 \\
        \vdots & \ddots & & \vdots\\
        \hdots & 0 & \mu & 0 \\
        c_{k_0} & \hdots & c_1 & \lambda
    \end{bmatrix}
    }_{Lf}
    \cdot
    \underbrace{
    \begin{bmatrix}
        y^{k_0}\\
        \vdots\\
        y\\
        x
    \end{bmatrix}}_{Lq}
    = L(f(q))
\end{align*}
The strictly subresonant ones are those for which $\lambda=\mu=1$.

We also illustrate the notion of cocycle normal forms.
Suppose that $E$ is a $2$-dimensional cocycle, smooth along stables, and with Lyapunov exponents $\mu_1>\mu_2$ and basis of linear functions $\xi_1,\xi_2$ (so, with weights $-\mu_1<-\mu_2$).
Then by a change of coordinates similar to that in \cref{thm:Normal-Forms-along-Stables} which is linear on the fibers, we can ensure that the cocycle takes the form
\begin{align*}
    \everymath{\displaystyle}
    \begin{bmatrix}
        \xi_1\\
        \xi_2\\
        x\\
        y
    \end{bmatrix}
    \mapsto
    \begin{bmatrix}
        a_1 \xi_1 + \xi_2\cdot\left(\sum_{i\lambda_1 + j\lambda_2\leq \mu_1-\mu_2} b_{i,j}x^i y^j\right)\\
        a_2 \xi_2 \phantom{+ \xi_2\cdot\left(\sum_{i\lambda_1 + j\lambda_2\leq \mu_1-\mu_2} b_{i,j}x^i y^j\right)}\\
        \lambda x + \sum_{0<k\lambda_2 \leq \lambda_1}c_k y^k \phantom{XXXXXXXX}\\
        \mu y \phantom{+ \sum_{0<k\lambda_2 \leq \lambda_1}c_k y^k}\phantom{XXXXXXXX}
    \end{bmatrix}
\end{align*}
Our convention is that the dynamics in normal form coordinates is centered at the origin.

%%%%% End of SubSection Classical Normal Forms %%%%%
%%%%%%%%%%%%%%%%%%%%%%%%%%%%%%%%%%%%%%%%%%%%%%%%%%%%%%%%%%%%%%%%%%%%%%%%%%%%%%%

%%%%%%%%%%%%%%%%%%%%%%%%%%%%%%%%%%%%%%%%%%%%%%%%%%%%%%%%%%%%%%%%%%%%%%%%%%%%%%%
%%%%% SubSection Normal Forms as (G,X)-structures %%%%%

\subsection{Normal Forms as \texorpdfstring{$(G,X)$}{(G,X)}-structures.}
        \label{ssec:Normal-Forms-as-G-X-structures}

Recall that given a group $G$ acting transitively on a ``model'' manifold $X$, one can define a $(G,X)$ structure on any manifold $M$ to consist of charts valued in open subsets of $X$, with gluing maps given by elements of $G$.
One can loosen the definition a bit further (without however increasing the class of spaces) by using several isomorphic model spaces.
With this in mind, we have:

\begin{definition}[Subresonant structure on a manifold]
	\label{def:subresonant_structure_on_a_manifold}
	Fix real scalars $-\lambda_1<\dots<-\lambda_n<0$ and appropriate multiplicities, denoted by $\overrightarrow{\lambda}$.
	A \emph{subresonant structure} of type $\overrightarrow{\lambda}$ on a manifold $M$, consists of charts $U_i\subseteq M$ to open subsets of filtered vector spaces $V_i$, where the filtration has type $\overrightarrow{\lambda}$, such that the transition maps are subresonant polynomial isomorphisms.

    Equivalently, we have the data of a $(G,X)$-structure with $X:=V$ a filtered vector space and $G:=\bbG^{sr}(V)$.
\end{definition}
A classical example of subresonant structure is an affine structure, namely manifolds with charts glued by the affine group $\GL_n(\bR)\ltimes \bR^n$.
This arises by taking only one scalar $-\lambda_1<0$ and $V^{\leq -\lambda_1}:=V$ is the only element of the filtration.

With this definition, we have a transparent formulation of \autoref{thm:Normal-Forms-along-Stables}: stable manifolds admit a subresonant structure, on which the dynamics is given by subresonant polynomial maps.
Let us see how this also gives an alternative approach to the existence of normal form coordinates, i.e. subresonant structures.

\noindent\textbf{Cartan Geometries.}
We first recall that a ``Cartan Geometry'' generalizes a $(G,X)$-structure in the same manner as Riemannian geometry generalizes Euclidean geometry, see \cite{Sharpe1997_Differential-geometry} for a general introduction.
Pick a point stabilizer $H\subset G$, so that $X\isom G/H$.
To give a Cartan geometry on a manifold $M$ is equivalent to giving a principal $H$-bundle $P\to M$, together with a $\frakg:=\Lie(G)$-valued $1$-form $\omega$ on $P$, subject to certain axioms (see \cite[E.4.2]{BrownEskinFilip2025_Measure-rigidity-for-generalized-u-Gibbs-states-and-stationary}).
A Cartan geometry has an associated notion of curvature, and the curvature vanishes if and only if the Cartan geometry is locally isomorphic to the principal $H$-bundle $G\to G/H\isom X$, i.e. a $(G,X)$-structure.

\noindent\textbf{Natural connections in dynamics.}
A basic principle in dynamics is the following: if $E\to Q$ is a smooth cocycle, $S\subset E$ is a subcocycle, and the Lyapunov exponents of $E/S$ are all strictly smaller than those of $S$, then we have an equivariant splitting $E/S\to E$.
Furthermore, along stable manifolds such a splitting is as regular as $S$ and $E$ (e.g. smooth if $S$ and $E$ are smooth).
Indeed, the stable Lyapunov filtration of $E$ is of the same regularity as $E$ along stable manifolds, and the assumption is that $S$ is equivariantly isomorphic to a corresponding stable Lyapunov subbundle of $E$.

One can apply this principle systematically, as is done in \cite[E.3-4]{BrownEskinFilip2025_Measure-rigidity-for-generalized-u-Gibbs-states-and-stationary}, to obtain two things:
\begin{itemize}
    \item At $\nu$-a.e. point $q\in Q$, the jet of the stable manifold, as well as the associated infinitesimal subresonant structure.
    \item A Cartan geometry modeled on $(\bbG^{sr}(V),V)$, where the connection $1$-form $\omega$ is determined by dynamical splittings as just described.
\end{itemize}
One then sees that the Cartan geometry has zero curvature, as the curvature is equivariant for the dynamics, which is contracting along the stable directions, and hence the curvature must vanish.
This yields the required dynamically equivariant subresonant structure on stable manifolds.

%%%%% End of SubSection Normal Forms as (G,X)-structures %%%%%
%%%%%%%%%%%%%%%%%%%%%%%%%%%%%%%%%%%%%%%%%%%%%%%%%%%%%%%%%%%%%%%%%%%%%%%%%%%%%%%%

%%%%%%%%%%%%%%%%%%%%%%%%%%%%%%%%%%%%%%%%%%%%%%%%%%%%%%%%%%%%%%%%%%%%%%%%%%%%%%%%
%%%%% SubSection Generalized u-Gibbs states %%%%%

    \subsection{Generalized u-Gibbs states.}
        \label{ssec:Generalized-u-Gibbs-states}

With the language of normal forms, it is now possible to give some of the key definitions of the theory.
We omit several technical points in the definition below and refer to \cite[Def.~3.1.5]{BrownEskinFilip2025_Measure-rigidity-for-generalized-u-Gibbs-states-and-stationary} for details.
For a point $q\in Q$ and subgroup $U^+(q)\subseteq \bbG^{ssr}(\cW^u[q])$ we denote by $U^+[q]:=U^+(q)\cdot q\subset \cW^u[q]$ the associated orbit.

\begin{definition}[Compatible family of subgroups]
    \label{def:compatible_family_of_subgroups}
    Suppose that $\nu$ is an ergodic probability measure on $Q$, invariant under a diffeomorphism $f\in \Diff(Q)$.
    A \emph{compatible family of subgroups} is the data of an $f$-equivariant measurable map $q\mapsto U^+(q)$ where $U^+(q)\subseteq \bbG^{ssr}(q)$ is an algebraic subgroup, with the following properties:
    \begin{enumerate}
        \item The partition with atoms $q\mapsto \frakB_0[q]\cap U^+[q]$ is $\nu$-measurable.
        \item The conditional measure on $U^+[q]$ obtained from this partition is $U^+(q)$-invariant, and $U^+(q)$ is the largest connected subgroup of $\bbG^{ssr}(\cW^u[q])$ that preserves it.
        \item For $q'\in U^+[q]$ we have that $U^+(q)=U^+(q')$ as subgroups of $\bbG^{ssr}(\cW^u[q]) = \bbG^{ssr}(\cW^u[q'])$.
    \end{enumerate}
\end{definition}
Informally, this definition allows us to speak of $f$-invariant measures that have Lebesgue-class conditionals along some nontrivial directions, not necessarily coming from the filtration by unstable submanifolds.
We denote by $\cB_0[q]:=\frakB_0[q] \cap U^+[q]$ the measurable partition from the definition.

\begin{definition}[Generalized u-Gibbs state]
    \label{def:Generalized-u-Gibbs-state}
    An ergodic $f$-invariant probability measure $\nu$ on $Q$ is a \emph{generalized u-Gibbs state} if there exists a nontrivial measurable family of subgroups $U^+(q)\subseteq \bbG^{ssr}(q)$ compatible with $\nu$, in the sense of \autoref{def:compatible_family_of_subgroups}.
\end{definition}

Recall (see Pesin--Sinai \cite{PesinSinai1982_Gibbs-measures-for-partially-hyperbolic-attractors}) that ``u-Gibbs states'' are measures whose conditionals along the fast unstables are of Lebesgue class, which are a special case of \cref{def:Generalized-u-Gibbs-state}.
It is always possible to construct u-Gibbs states, and it is important to know when they are unique, see e.g. Alvarez--Leguil--Obata \cite{AlvarezLeguilObata2022_Rigidity-of-U-Gibbs-measures-near-conservative-Anosov}.

\noindent\textbf{Linearization cocycle.}
Suppose given a compatible family of subgroups $U^+$.
Then one can introduce the space $\cC(q):=\bbG^{ssr}(q)/U^+(q)$ of translates of $U^+[q]$ inside $\cW^u[q]$.
It turns out that this space admits a subresonant structure and hence we have a linearization cocycle $L\cC$, with a dynamics-equivariant embedding $\cC\into L\cC$ (see \cite[\S5]{BrownEskinFilip2025_Measure-rigidity-for-generalized-u-Gibbs-states-and-stationary}).
The cocycle $L\cC$ has a tautological equivariant section given by the identity coset, and we denote by $\bbH$ the quotient.
This cocycle measures the divergence of nearby $U^+$-orbits.

\noindent\textbf{Finite cover.}
There is an important technical point which we have omitted in \cref{def:compatible_family_of_subgroups}.
Namely, the invariant family of subgroups might be defined on a finite cover of $(Q,\nu)$.
This arises because one might need to pass to a finite cover to bring the cocycle $\bbH$ to ``Jordan normal form''.
This has a ripple effect on most concepts below, but for notational simplicity we continue to work with the space $(Q,\nu)$ instead of its possible finite cover $X$ as in \cite{BrownEskinFilip2025_Measure-rigidity-for-generalized-u-Gibbs-states-and-stationary}.

%%%%% End of SubSection Generalized u-Gibbs states %%%%%
%%%%%%%%%%%%%%%%%%%%%%%%%%%%%%%%%%%%%%%%%%%%%%%%%%%%%%%%%%%%%%%%%%%%%%%%%%%%%%%%

%%%%%%%%%%%%%%%%%%%%%%%%%%%%%%%%%%%%%%%%%%%%%%%%%%%%%%%%%%%%%%%%%%%%%%%%%%%%%%%%
%%%%%%%%%%%% End of Section Normal Forms %%%%%%%%%%%%%%%%%%%%
%%%%%%%%%%%%%%%%%%%%%%%%%%%%%%%%%%%%%%%%%%%%%%%%%%%%%%%%%%%%%%%%%%%%%%%%%%%%%%%%

%%%%%%%%%%%%%%%%%%%%%%%%%%%%%%%%%%%%%%%%%%%%%%%%%%%%%%%%%%%%%%%%%%%%%%%%%%%%%%%%
%%%%%%%%%%%%%%%%%%%% Section Extra Invariance or Joint Integrability %%%%%%%%%%%%%%%%%%%%
%%%%%%%%%%%%%%%%%%%%%%%%%%%%%%%%%%%%%%%%%%%%%%%%%%%%%%%%%%%%%%%%%%%%%%%%%%%%%%%%

                    \section{Extra Invariance or Joint Integrability.}
                    \label{sec:Extra-Invariance-or-Joint-Integrability}

After introducing the crucial notion of ``Quantitative NonIntegrability'' in \cref{ssec:QNI-Quantitative-NonIntegrability}, we present the main theorem and an overview of the proof.

%%%%%%%%%%%%%%%%%%%%%%%%%%%%%%%%%%%%%%%%%%%%%%%%%%%%%%%%%%%%%%%%%%%%%%%%%%%%%%%%
%%%%% SubSection QNI (Quantitative NonIntegrability) %%%%%

    \subsection{QNI (Quantitative NonIntegrability).}
        \label{ssec:QNI-Quantitative-NonIntegrability}

We continue with notation as in the previous sections, but for convenience and compatibility of notation with \cite{BrownEskinFilip2025_Measure-rigidity-for-generalized-u-Gibbs-states-and-stationary}, we switch to continuous time.
So the dynamics is given by a flow $g_t$ instead of a map $f$.

Suppose that $\nu$ is a generalized u-Gibbs state, with associated measurable family of subgroups $U^+(q)\subseteq \bbG^{ssr}(q)$ compatible with $\nu$.
We can hope that $\nu$ might be of Lebesgue class along a submanifold of $Q$, but in general this can be obstructed by the fact that the measure could be, for example, the product of a Lebesgue measure along a submanifold and a fractal measure in the transverse direction.
To rule out this possibility, we introduce the following condition, which is the contrapositive of a joint integrability condition.
We have to be somewhat informal here, and refer to \cite[\S3.2]{BrownEskinFilip2025_Measure-rigidity-for-generalized-u-Gibbs-states-and-stationary} for a precise definition.

In the definition below, the notation for points $q_{\bullet}$ and time parameter $\ell$ is compatible with the outline of the proof in \cref{ssec:Scheme-of-Proof} and \cref{fig:Y_diagram} below.
Recall that $\frakB_0^-$ is a measurable partition subordinated to the stables, and $\cB_0$ is a measurable partition subordinated to the unstables and refined by the family $U^+$.

\begin{definition}[Quantitative NonIntegrability (QNI)]
    \label{def:Quantitative-NonIntegrability-QNI}
    We say that the compatible family of subgroups $U^+$ satisfies \emph{quantitative non-integrability} if there exist:
    \begin{itemize}
        \item a constant $\alpha_0>0$
        \item for any $\delta>0$ a compact set $K$ with $\nu(K)>1-\delta$, as well as constants $C,\ell_0>0$ depending on $\delta$
    \end{itemize}
    such that if $\ell>\ell_0$, and $q,q_{1/2},q_1\in K$ where $q_{1/2}:=g_{\ell/2}q$ and $q_1:=g_{\ell}q$, we have that
    \begin{itemize}
        \item there exists a measurable set $S\subset g_{-\ell/2}\left(\cB_0[q_1]\right)$ with $|S|/|g_{-\ell/2}\left(\cB_0[q_1]\right)|\geq 1-\delta$, where $|\bullet|$ denotes Haar measure on $U^+[q]$,
        \item there exists a measurable set $S'\subset g_{-\ell/2}\left(\frakB_0^-[q_1]\right)$ with $\nu^{s}{[q_{1/2}]}(S')/\nu^{s}{[q_{1/2}]}(g_{-\ell/2}\left(\frakB_0[q_1]\right))\geq 1-\delta$
    \end{itemize}
    such that for each $y_{1/2}\in S$ and $q_{1/2}'\in S'$ we have that
    \[
        d\left(\cB_0[q'_{1/2}],\cW^{cs}_{loc}[y_{1/2}]\right)\geq C e^{-\alpha_0 \ell}.
    \]
\end{definition}
We note that we work with an explicit realization of the center-stable manifold above, and the exponent $\alpha_0$ is not affected by this choice.
Implicit in the definition is that $U^+$ is defined at all points of $S'$ and $\cW^{cs}_{loc}$ at all points of $S$.
The condition on the measure of $S$ and $S'$ (which is a density requirement) could have been equivalently stated at $q_1$ to be $|g_{\ell/2}S|/|\cB_0[q]|\geq 1-\delta$ and similarly for $S'$.
See also Katz \cite{Katz2023_Measure-rigidity-of-Anosov-flows-via-the-factorization-method} and Eskin--Potrie--Zhang \cite{EskinPotrieZhang2023_Geometric-properties-of-partially-hyperbolic-measures-and-applications} for related works involving the QNI condition.

%%%%% End of SubSection QNI (Quantitative NonIntegrability) %%%%%
%%%%%%%%%%%%%%%%%%%%%%%%%%%%%%%%%%%%%%%%%%%%%%%%%%%%%%%%%%%%%%%%%%%%%%%%%%%%%%%%

%%%%%%%%%%%%%%%%%%%%%%%%%%%%%%%%%%%%%%%%%%%%%%%%%%%%%%%%%%%%%%%%%%%%%%%%%%%%%%%%
%%%%% SubSection Main Theorem %%%%%

    \subsection{Main Theorem.}
        \label{ssec:Main-Theorem}

With all these preliminaries, we are now ready to state our main result:

\begin{theorem}[Inductive step]
    \label{thm:Inductive-step}
    Suppose that $\nu$ is an ergodic $f$-invariant probability measure on $Q$, which is a generalized u-Gibbs state with associated measurable family of subgroups $U^+(q)\subseteq \bbG^{ssr}(q)$ compatible with $\nu$.
    If $\nu$ satisfies the QNI condition (\cref{def:Quantitative-NonIntegrability-QNI}), then there exists a strictly larger-dimensional measurable family of subgroups $\widetilde{U}^+(q)\subseteq \bbG^{ssr}(q)$ compatible with $\nu$.
\end{theorem}
In other words, if $U^+(q)$ is as large as possible, then $\nu$ does not satisfy the QNI condition.
We can then expect to be able to integrate the stable and unstable directions and show that $\nu$ is of Lebesgue class along such a foliation.

%%%%% End of SubSection Main Theorem %%%%%
%%%%%%%%%%%%%%%%%%%%%%%%%%%%%%%%%%%%%%%%%%%%%%%%%%%%%%%%%%%%%%%%%%%%%%%%%%%%%%%%

%%%%%%%%%%%%%%%%%%%%%%%%%%%%%%%%%%%%%%%%%%%%%%%%%%%%%%%%%%%%%%%%%%%%%%%%%%%%%%%%
%%%%% SubSection Scheme of Proof %%%%%

    \subsection{Scheme of Proof.}
        \label{ssec:Scheme-of-Proof}

We include an overview of some of the main steps that go into the proof of \cref{thm:Inductive-step}.
Suppose that we are given a nontrivial compatible family of subgroups $U^+(q)$ satisfying QNI.
Recall that $\nu^u[q]$ denotes the conditional measure along unstables.
We would like to show that, for any $\ve>0$, there exists a positive measure set of points $x,y$, with $y\in \cW^u[x]$ and such that:
\begin{enumerate}
  \item We have separation of $U^+$-orbits, where $d_{\cH,loc}$ denotes the (local) Hausdorff distance between two sets:
  \[
    d_{\cH,loc}(U^+[x],U^+[y]) \asymp \ve.
  \]
  \item We have invariance of conditional measures (up to scaling, where $\propto$ denotes two measures that agree up to a constant factor):
  \begin{align}
    \label{eqn:main_invariance_property}
    \nu^u[y] \propto \wp^u(x,y)\nu^u[x]
  \end{align}
  
  where $\wp^u(x,y)\in \bbG^{ssr}(\cW^u[x])$ is a strictly subresonant map given by the unstable measurable connection discussed in \cref{ssec:Classical-Normal-Forms}.
\end{enumerate}
As $\ve\to 0$, the last condition would yield invariance of $\nu^u$ by a $1$-parameter group of strictly subresonant transformations, while the first condition implies that this subgroup is not contained in $U^+$.

The actual argument shows invariance of a finer family of conditional measures (see \cite[Prop.~12.0.2]{BrownEskinFilip2025_Measure-rigidity-for-generalized-u-Gibbs-states-and-stationary}) denoted $\wtilde{f}_{ij}$, which will be explained below.

\noindent\textbf{Some conventions and notation.}
Throughout, one frequently passes to a compact set $K$ of almost full measure, on which the ergodic theorems hold uniformly for the relevant observables, the measurably varying objects are uniformly continuous, and the involved constants are uniformly bounded.
We will not specify this repeated passage to smaller compact sets and always denote by $K$ the set which is relevant for an argument.

A quantity $X$ (or vector, or distance) will be called ``exponentially small'' if there exists $\delta_X>0$ depending only on the Lyapunov spectrum and the QNI condition, and a constant $C$ (depending in addition on the compact set $K$) such that $|X|\leq Ce^{-\delta_X \ell}$, where $\ell$ is the main time parameter introduced below.
A quantity will be called ``macroscopic'' if $\tfrac 1 C\leq |X|\leq C$ again for a constant $C$ depending on the compact set $K$.
Both expressions can be abbreviated as $|X|\asymp_K 1$ and $|X|\leqapprox_K e^{-\delta_X \ell}$.
More generally $X\leqapprox_S Y$ will mean that there exists a constant $C$ depending on $S$ such that $X\leq C\cdot Y$.
We will write $X\leqapprox_\lambda Y$ when the implicit constant depends only on the Lyapunov spectrum and exponent in the QNI condition.

\noindent\textbf{The eight points.}
The scheme for obtaining \cref{eqn:main_invariance_property} goes back to the work of Benoist--Quint and Eskin--Mirzakhani \cite{BenoistQuint2011_Mesures-stationnaires-et-fermes-invariants-des-espaces-homogenes,EskinMirzakhani2018_Invariant-and-stationary-measures-for-the-rm-SL2R-action-on-moduli-space}.
% To be consistent with \cite{BrownEskinFilip2025_Measure-rigidity-for-generalized-u-Gibbs-states-and-stationary}, we will work with continuous time and denote the flow by $g_t$.
To start, suppose we have $q,q'\in K$ with $q'\in \cW^s_{loc}[q]$ at a macroscopic distance. 
A crucial time parameter will be used throughout the argument and denoted by $\ell$, it will lead to two additional time parameters $t,\tau$, each satisfying $t\leqapprox_{\lambda}\ell$ and $\tau\leqapprox_{\lambda}\ell$.

Introduce next the points
\begin{align*}
    \begin{split}
    q_1 & :=g_\ell q\\
    q_1' & :=g_{\ell}q_1
    \end{split}
    \begin{split}
    q_3  & :=g_t q_1\\
    q_3' & :=g_t q_1'
    \end{split}
\end{align*}
Note that $q_1$ and $q_1'$ are exponentially close (since they are stably related) and $q_3,q_3'$ are even closer.
We continue to assume that all the points are in $K$.

\begin{figure}[htbp!]
    \centering
    \includegraphics[width=0.4\linewidth]{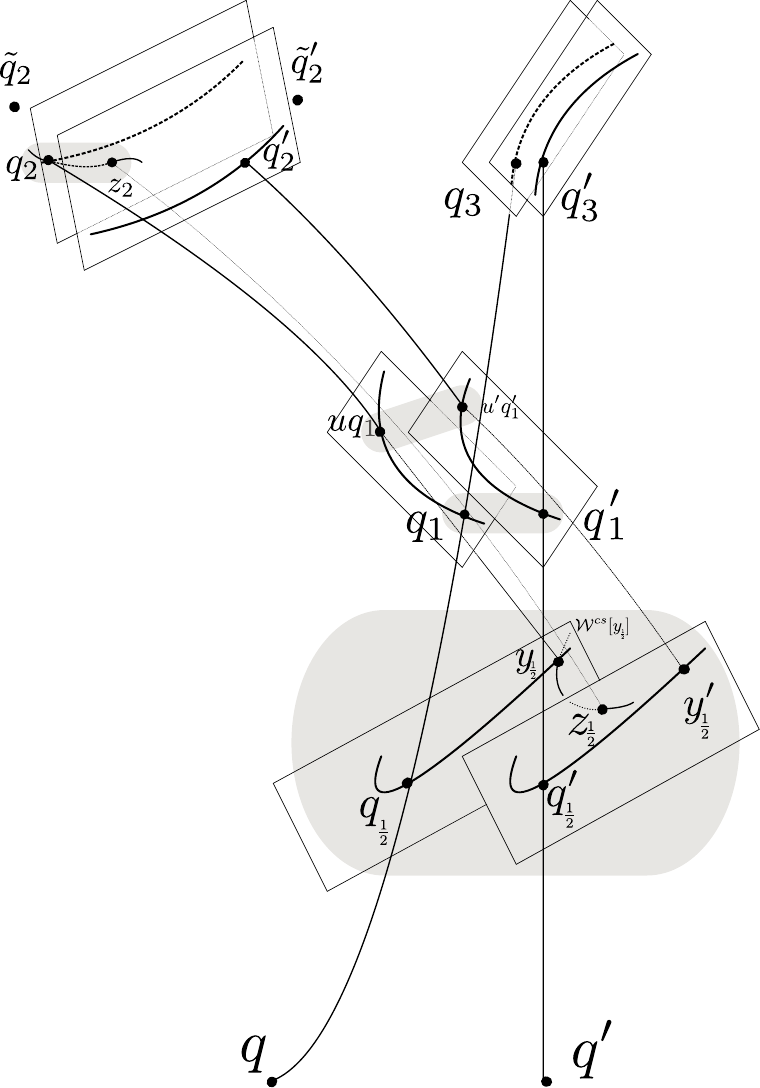}
    \caption{Configuration of points in the main argument.
    Points in the same shaded cloud are exponentially close.}
    \label{fig:Y_diagram}
\end{figure}

Next, choose $u\in U^+(q_1),u'\in U^+(q_1')$, of macroscopic size but close to each other (after appropriate identifications of $U^+$-groups) and set
\[
    q_2:=g_{\tau}(u q_1) \quad q_2' := g_{\tau}(u'q_1')
\]
where $uq_1$ denotes the result of acting by $u\in U^+(q_1)\subset \bbG^{ssr}(\cW^u[q_1])$ on $q_1$, and analogously for the primed points, see \cref{fig:Y_diagram}.

\noindent\textbf{Exponential Drift.}
We let $\ell\to +\infty$ and pick $\tau$ such that $U^+[q_2]$ and $U^+[q_2']$ are (locally) at a Hausdorff distance $\asymp_K \ve$, and by passing to a subsequence can assume that $q_2\to \wtilde{q}_2$ and $q_2'\to \wtilde{q}_2'$ with $\wtilde{q}_2,\wtilde{q}_2'\in K$ and $\wtilde{q}_2'\in \cW^u[\wtilde{q}_2]$.

Recall that the conditional measures are equivariant, and the dynamics is by subresonant maps.
Therefore, we have that
\begin{align}
    \label{eqn:up_and_down_the_Y_diagram}
    \begin{split}
    \nu^u[q_2] & \propto A_* \nu^u[q_3]\\
    \nu^u[q_2'] & \propto A_*' \nu^u[q_3']
    \end{split}
\end{align}
for subresonant maps $A,A'$ between the corresponding unstables.
By picking the time parameters $t,\tau,\ell$ accordingly, we can ensure that $A$ and $A'$ are comparable and of bounded size (in appropriate trivializations).

In the limit as $\ell\to +\infty$ we have that $\nu^u[q_3]$ approaches $\nu^u[q_3']$ because $q_3\to q_3'$ (again, in appropriate trivializations) whereas the difference $A'\circ A^{-1}$ converges to the strictly subresonant map $\wp^u(\wtilde{q}_2,\wtilde{q}_2')$ so we obtain
\[
    \nu^u[\wtilde{q}_2'] \propto \wp^u(\wtilde{q}_2,\wtilde{q}_2')\nu^u[\wtilde{q}_2]
\]
as desired.

\noindent\textbf{Divergence and conditionals.}
For the above argument to yield nontrivial invariance, it is important that $uq_1$ and $u'q_1'$ diverge under the forward time dynamics (in particular, are not on the same stable manifold).
In fact, we need divergence of the $U^+$-orbits.
We recall from \cref{ssec:Generalized-u-Gibbs-states} that we have a linear cocycle $\bbH$ that controls the divergence of $U^+$-orbits.
It is then possible to construct a nonempty subcocycle $\bbE\subset \bbH$ (see \cite[\S10.1]{BrownEskinFilip2025_Measure-rigidity-for-generalized-u-Gibbs-states-and-stationary}) such that under the process of going up and down the top half of \cref{fig:Y_diagram}, for most choices of $u$, the divergence is along $\bbE$.

The Oseledets decomposition of $\bbE$ is then further refined to arrange the cocycle dynamics to be block-conformal.
This is then used to further refine the partition by unstable manifolds and disintegrate the measures $\nu^u$ to a finer family $\wtilde{f}_{ij}$, for which actual extra invariance is obtained (see \cite[\S10-11]{BrownEskinFilip2025_Measure-rigidity-for-generalized-u-Gibbs-states-and-stationary}).
Note that it is this step that forces the finite cover mentioned in \cref{ssec:Generalized-u-Gibbs-states}.

\noindent\textbf{Factorization.}
In the above argument, it is important to have accurate estimates of various objects, such as the orbits $U^+[q]$ and their divergence, as well as accurate comparisons of these objects at nearby points.
A crucial tool to accomplish this is the introduction of the ``half-way points'' $q_{1/2}:=g_{\ell/2}q,q_{1/2}':=g_{\ell/2}q'$ as in \cref{fig:Y_diagram}.
The reason these are useful is that at that region of the diagram, all the points that appear at later stages in the argument are exponentially close.
This means that the Taylor expansion of various quantities becomes useful, and higher orders are even more useful as we can approximate various quantities to arbitrarily good exponential accuracy.

A formalism for doing such approximations is developed in \cite[\S7]{BrownEskinFilip2025_Measure-rigidity-for-generalized-u-Gibbs-states-and-stationary}, and quantities that can be well-approximated, in a manner compatible with the dynamics, are called ``factorizable''.

\noindent\textbf{Invariance and Entropy.}
We end with a sketch of proof of \cref{thm:stiffness}, with notation as in \cref{ssec:Setup-Manifolds-Flows-Random-Walks} and \cref{ssec:Invariant-Manifolds}.
Applying iteratively \cref{thm:Inductive-step} to a $\mu$-stationary measure $\nu$ yields a maximal family of groups $U^+$ compatible with the measure.
We can define $\cL^s[\hat{q}]\subseteq \cW^s[\hat{q}]$ to be the Zariski-closure of $\supp \nu^s[\hat{q}]$ and the bundle $\cZ(q)\subset T_qQ$ (for $\nu$-a.e. $q$) as the span of $T_q \cL^s[q,b]$ over $\mu$-a.e. choice of $b\in \Diff(Q)^\bZ$.
Clearly $\cZ$ is invariant under the random walk, so the sum of its Lyapunov exponents is nonnegative.

Furthermore, at this stage we also know that $\nu^u[q,b]$ is of Lebesgue class in the directions $\cZ(q)\cap T_q U^+[q,b]$ so we find from the Ledrappier--Young formula \cref{eqn:Ledrappier_Young_Variant} that
\[
    h_{\hat{\nu}}\left(\hat{F}\right)
    \geq
    H(\mu) + 
    \sum_{\lambda_i>0} \lambda_i\cdot \dim \left(\cZ^{\lambda_i}\right)
\]
where $H(\mu)$ is the entropy of the probability distribution on the support of $\mu$.
Reversing time direction, and using that $\supp \nu^s[q]\subset \cL^s[q]$, we find:
\begin{align*}
    h_{\hat{\nu}}\left(\hat{F}^{-1}\right) 
    & \leq
    h_{\hat{\nu}}\left(\hat{F}^{-1};\frakW^s_1\right)
    + \sum_{\lambda_i<0} (-\lambda_i) \cdot \dim \left(\cL^{s,\lambda_i}\right)\\
    & \leq
    h_{\hat{\nu}}\left(\hat{F}^{-1};\frakW^s_1\right)
    + \sum_{\lambda_i<0} (-\lambda_i) \cdot \dim \left(\cZ^{\lambda_i}\right)\\
\end{align*}

Using that $h_{\hat{\nu}} \left(\hat{F}\right) =h_{\hat{\nu}}\left(\hat{F}^{-1}\right)$ we find that
\[
    h_{\hat{\nu}}(\hat{F}^{-1};\frakW^s_{1}) - H(\mu)
    \geq
    \sum_{\lambda_i}\lambda_i \cdot \dim \left(\cZ^{\lambda_i}\right)\geq 0.
\]
Using that $H(\mu)\geq h_{\hat{\nu}}\left(\hat{F}^{-1};\frakW^{s}_{1}\right)$ if and only if $\nu$ is $\mu$-invariant (the Ledrappier invariance principle \cite[Prop.~2]{Ledrappier1986_Positivity-of-the-exponent-for-stationary-sequences-of-matrices}, see also \cref{eqn:relative_entropy_bound}) the result follows.

%%%%% End of SubSection Scheme of Proof %%%%%
%%%%%%%%%%%%%%%%%%%%%%%%%%%%%%%%%%%%%%%%%%%%%%%%%%%%%%%%%%%%%%%%%%%%%%%%%%%%%%%

%%%%%%%%%%%%%%%%%%%%%%%%%%%%%%%%%%%%%%%%%%%%%%%%%%%%%%%%%%%%%%%%%%%%%%%%%%%%%%%
%%%%%%%%%%%% End of Section Extra Invariance or Joint Integrability %%%%%%%%%%%%%%%%%%%%
%%%%%%%%%%%%%%%%%%%%%%%%%%%%%%%%%%%%%%%%%%%%%%%%%%%%%%%%%%%%%%%%%%%%%%%%%%%%%%%

%%%%%%%%%%%%%%%%%%%%%%%%%%%%%%%%%%%%%%%%%%%%%%%%%%%%%%%%%%%%%%%%%%%%%%%%%%%%%%%
%%%%%%%%%%%%%%%%%%%% Section Further Directions %%%%%%%%%%%%%%%%%%%%
%%%%%%%%%%%%%%%%%%%%%%%%%%%%%%%%%%%%%%%%%%%%%%%%%%%%%%%%%%%%%%%%%%%%%%%%%%%%%%%

                    \section{Further Directions.}
                    \label{sec:Further-Directions}

We conclude with a series of questions that naturally result from the results we've discussed.

%%=============================================================================
%%                    start of subsec: Questions about orbit closures and measures

\subsection{Questions about orbit closures and measures.}
    \label{ssec:questions_about_orbit_closures_and_measures}
\hfill

\noindent\textbf{Empirical measures.}
We remain in the random walk setting of \cref{ssec:Setup-Manifolds-Flows-Random-Walks}.
Denote by $\Gamma_\mu$ the subgroup of $\Diff(Q)$ generated by the support of $\mu$.
Given a point $q\in Q$, one can consider the empirical measures
\[    
    \nu_{N,q} := \frac 1 N \sum_{n=0}^{N-1} \mu^{*n} * \delta_q
\]
Any weak-$*$ limit of $\nu_N$ is a $\mu$-stationary measure, and in favorable situations one can hope that the support of the limit measure coincides with the orbit closure $\overline{\Gamma_\mu\cdot q}$ (it is always contained in the orbit closure by construction).
We will refer to a stationary measure arising this way as ``associated to $q$''.
The typical approach to classifying orbit closures is to understand first the ergodic invariant measures.

\noindent\textbf{Spaces stratified by manifolds.}
In settings where one expects a classification of invariant measures and orbit closures, these should admit an inductive structure where the main stratum is a smooth manifold, and its closure consists of orbit closures of lower dimension.
At a minimum, we hope for:

\begin{definition}[Space stratified by manifolds]
    \label{def:space_stratified_by_manifolds}
    A subset $N$ in a manifold $Q$ will be said to be \emph{stratified by manifolds} if it is endowed with a finite decomposition $N=\coprod S_i$, with each $S_i$ connected and called a \emph{stratum}, satisfying the following properties:
    \begin{enumerate}
        \item For each $S_i$ and point $s\in S_i$, there exists an open $U\ni s$ such that $S_i\cap U$ is an embedded submanifold of $U$.
        \item If the closure of a stratum $\ov{S_i}$ intersects another stratum $S_j$, then $S_j\subset \ov{S_i}$ and $\dim S_j<\dim S_i$.
    \end{enumerate}
    Denote by $N^{\circ}\subset N$ the union of strata of top dimension, and call this top dimension the \emph{dimension of $N$}.
\end{definition}
There is a plethora of further constraints that one could impose on the local structure and how strata are ``glued'' near singularities, but we omit this discussion as different settings might involve different further properties.
We refer to Mather's notes \cite{MR2958928} and the foreword by Goresky for an overview of the many issues involved, but would like to point out a common theme with that in dynamical systems: one wants to identify an open and dense set of mappings and singularities which admit a good structure theory.

\noindent\textbf{Smooth vs. Real Analytic.}
It is worth pointing out that the questions we pose in this section are likely to be more easily approached in the case of real analytic dynamical systems, rather than just smooth.
For example, it would be very interesting to relax the analyticity condition in \cref{thm:Totally-Geodesic-Rigidity} on totally geodesic hypersurfaces.
In the real-analytic case, one would expect that a ``space stratified by manifolds'' is at least a Whitney stratification.

\noindent\textbf{The case of K3 surfaces.}
A rich setting for many questions in dynamics is provided by K3 surfaces: a class of compact complex $2$-dimensional manifolds that admit a canonical (up to scale) holomorphic complex symplectic form.
In particular, these admit large automorphism groups, see \cite{Filip2022_An-introduction-to-K3-surfaces-and-their-dynamics} for an introduction to the topic and more context.

Suppose that $X$ is a complex K3 surface and let $\Aut(X)$ be its automorphism group.
A subgroup of $\Aut(X)$ is called \emph{non-elementary} if it contains a nonabelian free group.
A real $2$-dimensional submanifold $S\subset X$ is called \emph{totally real} if for every $s\in S$ the complex span of the real tangent space $T_sS$ agrees with $T_s X$.
The work of Cantat--Dujardin \cite{CantatDujardin2023_Random-dynamics-on-real-and-complex-projective-surfaces,CantatDujardin2023_Dynamics-of-automorphism-groups-of-projective-surfaces:-classification-examples,CantatDujardin2024_Finite-orbits-for-large-groups-of-automorphisms-of-projective} and Roda \cite{Roda2024_Classifying-hyperbolic-ergodic-stationary-measures-on-compact} suggests the following:

\begin{conjecture}[Orbit of a point on a K3 surface]
     \label{cjc:orbit_of_a_point_on_a_k3_surface}
    Assume that $\Gamma\subseteq \Aut(X)$ is a nonelementary subgroup and $\mu$ is a probability measure on a generating set of $\Gamma$.
    For any $x\in X$, exactly one of the following possibilities occurs:
    \begin{enumerate}
        \item The orbit $\Gamma\cdot x$ is finite.
        % , and the associated stationary measure is atomic.
        \item The orbit $\Gamma\cdot x$ is contained and Zariski-dense in an algebraic curve $C$.
        \item There exists $\Gamma$-invariant set $S\subset X$ containing $x$, such that $S$ is stratified by manifolds, is $2$-dimensional, the open stratum $S^{\circ}$ is totally real, and $\Gamma\cdot x$ is dense in $S$.
        The associated stationary measure is a smooth volume on $S^{\circ}$.
        \item The orbit $\Gamma\cdot x$ is dense in $X$, and the associated stationary measure is the smooth volume on it.
    \end{enumerate}
    Any ergodic stationary measure is obtained as an empirical measures from a point as above, and in all cases except that of an invariant algebraic curve the stationary measure is also invariant.
\end{conjecture}
In the first case above, the associated measure is atomic, and in the second case the curve $C$ is isomorphic to a finite union of $\bP^1(\bC)$'s, the action of $\Gamma$ on a component of $C$ is by a subgroup of $\PGL_2(\bC)$, and the associated stationary measure is determined by this linear action using Furstenberg's theory.
Note also it is possible for $x$ to be contained in a totally real surface but not have dense orbit, in which case the above conjecture says that the orbit is either finite, or confined to a real algebraic curve.

When $\Gamma_\mu$ contains parabolic elements, \cref{cjc:orbit_of_a_point_on_a_k3_surface} has been established by Cantat--Dujardin \cite[Thm.~1.22]{CantatDujardin2023_Dynamics-of-automorphism-groups-of-projective-surfaces:-classification-examples}.
The tools necessary to address the case of a totally real surface, in full generality, have been developed by Roda \cite{Roda2024_Classifying-hyperbolic-ergodic-stationary-measures-on-compact}.
The following special case would be a substantial progress towards Conjecture~\ref{cjc:orbit_of_a_point_on_a_k3_surface}, and most likely require substantial new ideas:

\begin{conjecture}[Measure with zero Lyapunov exponents on K3]
	\label{cjc:measure_with_zero_lyapunov_exponents_on_k3}
	With assumptions as in Conjecture~\ref{cjc:orbit_of_a_point_on_a_k3_surface}, suppose $\nu$ is a $\mu$-stationary measure with zero Lyapunov exponents.

	Then the support of $\nu$ is not Zariski-dense in $X$.
\end{conjecture}
This is known by work of Cantat--Dujardin when $\Gamma_\mu$ contains parabolics, and the difficulty is to remove this assumption.
Note that a $\mu$-stationary measure with zero Lyapunov exponents must necessarily be $\Gamma_\mu$-invariant, by Ledrappier's Invariance Principle \cite{Ledrappier1986_Positivity-of-the-exponent-for-stationary-sequences-of-matrices}, as extended by Crauel \cite{Crauel1993_Non-Markovian-invariant-measures-are-hyperbolic}.
Not having Zariski-dense support implies that the measure is either atomic, or supported on a finite union of $\bP^1(\bC)$'s.

\noindent\textbf{Symplectic Manifolds.}
The following question appearing in \cite[Conj.~1.1.12]{BrownEskinFilip2025_Measure-rigidity-for-generalized-u-Gibbs-states-and-stationary} is a natural extension of \cref{thm:Finite-or-Uniform}:

\begin{conjecture}[Orbit closures on symplectic manifolds]
    \label{cjc:orbit_closures_on_symplectic_manifolds}
    Suppose that $Q$ is a compact symplectic manifold and $\mu$ is supported on a group $\Gamma_\mu$ inside the smooth symplectomorphisms of $Q$.
    Assume that $\mu$ has uniform expansion in the sense of \cref{def:Uniform-Expansion-and-Uniform-Gaps} on every isotropic subspace of $T_qQ$, for every $q\in Q$, as well as uniform gaps in dimension $\tfrac 12 \dim Q$.
    
    Then every $\mu$-stationary measure $\nu$ on $Q$ is supported on a space stratified by manifolds $S_{\mu}\subseteq Q$, such that $S_\mu^{\circ}$ is a symplectic submanifold and $\nu$ is the volume measure on $S_{\mu}^{\circ}$.
    Furthermore, for any $q\in Q$ the orbit closure $\overline{\Gamma_\mu\cdot q}$ is such an $S_q\subseteq Q$, and the empirical measures of $q$ tend to the smooth volume on $S_q^{\circ}$.
\end{conjecture}
Uniform gaps in dimension $\tfrac 12 \dim Q$ ensures that there are no zero Lyapunov exponents and stable/unstable manifolds have uniform size, so that one can use the method of Dolgopyat--Krikorian \cite[\S10]{DolgopyatKrikorian2007_On-simultaneous-linearization-of-diffeomorphisms-of-the-sphere} to establish ergodicity of the natural volume.
The assumption of uniform expansion on isotropic subspaces should eliminate examples such as the projective actions on invariant algebraic curves on K3 surfaces, c.f. \cref{cjc:orbit_of_a_point_on_a_k3_surface}.

%%                    end of subsec: Questions about orbit closures and measures
%%=============================================================================

%%%%%%%%%%%%%%%%%%%%%%%%%%%%%%%%%%%%%%%%%%%%%%%%%%%%%%%%%%%%%%%%%%%%%%%%%%%%%%%
%%%%% SubSection Hyperbolicity of Stationary Measures %%%%%

    \subsection{Hyperbolicity of Stationary Measures.}
        \label{ssec:Hyperbolicity-of-Stationary-Measures}

Essentially all the measure rigidity results mentioned in this text need, as a starting point, a positive Lyapunov exponent of the stationary measure, see e.g. 
\cite{BenoistQuint2011_Mesures-stationnaires-et-fermes-invariants-des-espaces-homogenes,BrownRodriguez-Hertz2017_Measure-rigidity-for-random-dynamics-on-surfaces-and-related-skew,CantatDujardin2023_Random-dynamics-on-real-and-complex-projective-surfaces,EskinLindenstrauss2018_Random-walks-on-locally-homogeneous-spaces}.
In the homogeneous setting, this is well-understood and addressed by the theory of random matrix products initiated by Furstenberg \cite{Furstenberg1963_Noncommuting-random-products}, with many subsequent contributions.

In the general case of dynamics on manifolds, this remains mysterious.

\noindent\textbf{Uniform Expansion.}
Recall that in \cref{def:Uniform-Expansion-and-Uniform-Gaps} we introduced the notion of uniform expansion.
Some criteria for verifying it have been developed by Chung \cite{Chung2020_Stationary-measures-and-orbit-closures-of-uniformly-expanding-random} and these can be used to verify it numerically (see also earlier work of Liu \cite{Liu2016_Lyapunov-Exponents-Approximation-Symplectic-Cocycle-Deformation}).
Let us note that uniform expansion is an open property, i.e. persists under sufficiently small perturbations of the diffeomorphisms in the support.
The prevalence of uniform expansion under perturbations of the measures has been obtained by Elliot Smith \cite{Elliott-Smith2023_Uniformly-expanding-random-walks-on-manifolds}.

Nonetheless, it is expected that uniform expansion should hold, unless an ``obvious'' obstruction precludes it.
We thus pose:
\begin{problem}[Criteria for Uniform Expansion]
    \label{pb:criteria_for_uniform_expansion}
    Develop criteria for uniform expansion of a measure $\mu$ on $\Diff(Q)$, akin to those obtained for random walks on linear groups, see e.g. \cite{Furman2002_Random-walks-on-groups-and-random-transformations}.
\end{problem}
Such criteria should involve structures that are invariant under the entire group $\Gamma_{\mu}$, see \cite[Thm.~1.21]{CantatDujardin2023_Dynamics-of-automorphism-groups-of-projective-surfaces:-classification-examples} for a good example.
While it is not hard to obtain such \emph{measurable} structures, it would be very useful to improve some of these to higher regularity.
Let us mention that in the presence of uniform expansion, DeWitt--Dolgopyat \cite{DeWittDolgopyat2025_Conservative-Coexpanding-on-Average-Diffeomorphisms} obtained exponential mixing of the stationary measures.

%%%%% End of SubSection Hyperbolicity of Stationary Measures %%%%%
%%%%%%%%%%%%%%%%%%%%%%%%%%%%%%%%%%%%%%%%%%%%%%%%%%%%%%%%%%%%%%%%%%%%%%%%%%%%%%%

%%=============================================================================
%%                    start of subsec: Rigidity and Special Orbit Closures

\subsection{Rigidity and Special Orbit Closures.}
    \label{ssec:rigidity_and_special_orbit_closures}

This section will be deliberately vague about some concepts but will illustrate them using some of the conjectures from \cref{ssec:questions_about_orbit_closures_and_measures}.
Suppose that we are in a situation where we have a classification of invariant measures, as well as orbit closures (and all stationary measures are invariant).

We assume for this section that a ``good'' orbit closure $N$ is a space stratified by manifolds, and that the associated invariant measure $\nu$ is a smooth volume form on the smooth stratum $N^{\circ}$ of maximal dimension.
The need to separate orbit closures into ``good'' and ``bad'' arises already in the setting of K3 surfaces, as in e.g. Conjecture~\ref{cjc:orbit_of_a_point_on_a_k3_surface}.

\begin{definition}[Maximal Orbit Closure and Invariant Measure]
    \label{def:maximal_orbit_closure_and_invariant_measure}
    Suppose $(N,\nu)$ is a ``good'' orbit closure as above.
    We call $N$, resp. $\nu$ \emph{maximal} if there does not exist another such pair $(N',\nu')$ with $N\subsetneq N'$ and $\dim N<\dim N'<\dim Q$.
\end{definition}

With this in mind, we can now formulate the basic question:
\begin{problem}[Infinitely many implies Homogeneous]
    \label{pb:infinitely_many_implies_homogeneous}
    Suppose that $Q$ admits infinitely many good maximal orbit closures $Q_i$.
    Show that $Q$ admits a homogeneous structure preserved by $\Gamma_\mu$, and the invariant measures associated to $Q_i$ equidistribute to the measure associated to $Q$.
\end{problem}

\begin{remark}
    \label{rmk:on_homogeneity}
    A few comments about \cref{pb:infinitely_many_implies_homogeneous}:
    \begin{enumerate}
        \item We present this problem as a principle, rather than a statement that is expected to always hold.
        \item It is often more useful to think of the principle as a \emph{finiteness} statement: $Q$ has finitely many good maximal orbit closures, unless it has a structure that naturally accounts for them.
        \item The conclusion about homogeneity can be interpreted to mean that $Q$ itself can be stratified by manifolds in a $\Gamma_\mu$-invariant way, such that $Q^{\circ}$ admits a $(G,X)$-structure preserved by $\Gamma_\mu$, or perhaps more generally a Gromov rigid geometric structure \cite{Gromov1988_Rigid-transformations-groups}.
    \end{enumerate}
\end{remark}

We now discuss some evidence towards an affirmative answer to \cref{pb:infinitely_many_implies_homogeneous}.
\vskip\baselineskip

\noindent\textbf{Rigidity for K3 surfaces.}
In the case of K3 surfaces, a special case of \cref{pb:infinitely_many_implies_homogeneous} is made precise and proved by Cantat--Dujardin \cite{CantatDujardin2024_Finite-orbits-for-large-groups-of-automorphisms-of-projective}.
The strategy, typical in arithmetic dynamics, should be applicable in other situations where $\Gamma_{\mu}$ preserves an algebraic structure.

Specifically, if $X$ is a projective algebraic surface over $\ov{\bQ}$, the subgroup $\Gamma_\mu\subseteq \Aut(X)$ is nonelementary and contains parabolic elements, and $\Gamma_\mu$ has infinitely many periodic orbits, then the surface must be a Kummer example, i.e. obtained from a torus by blowing up finitely many points and taking a finite quotient.
A key point is that a periodic orbit is invariant by the Galois group.
For a fixed loxodromic element $\gamma\in \Gamma_\mu$, increasing sets of periodic points that are invariant by the Galois group equidistribute towards the measure of maximal entropy for $\gamma$.
Therefore, the empirical measures of the finite orbits equidistribute towards on the one hand the invariant volume, and on the other to the measure of maximal entropy for any loxodromic $\gamma\in \Gamma_\mu$.
This yields the Kummer property by the rigidity theorem obtained by Cantat--Dupont \cite{CantatDupont2020_Automorphisms-of-surfaces:-Kummer-rigidity-and-measure-of-maximal-entropy} (see \cite{FilipTosatti2021_Kummer-rigidity-for-K3-surface-automorphisms-via-Ricci-flat-metrics} for an alternative proof of Kummer rigidity).

Let us note that as far as we know, the case of \cref{pb:infinitely_many_implies_homogeneous} when a K3 surface contains infinitely many totally real $\Gamma_\mu$-invariant surfaces, is not addressed in the literature, see however \cite[Thm.~1.22]{CantatDujardin2023_Dynamics-of-automorphism-groups-of-projective-surfaces:-classification-examples} and \cite[Thm.~10.1]{CantatDujardin2025_Hyperbolicity-for-large-automorphism-groups-of-projective-surfaces} for finiteness statements.
For algebraic group actions on higher-dimensional Calabi--Yau manifolds, it should be possible to extend the construction of canonical heights from \cite{FilipTosatti2023_Canonical-currents-and-heights-for-K3-surfaces} and use these heights to study finite orbits and invariant subvarieties.

\noindent\textbf{Totally Geodesic Submanifolds.}
Another result that provides some evidence towards \cref{pb:infinitely_many_implies_homogeneous} is \cref{thm:Totally-Geodesic-Rigidity} on totally geodesic hypersurfaces in a negatively curved Riemannian manifold.
In this setting, \cref{pb:infinitely_many_implies_homogeneous} predicts:

\begin{conjecture}[Totally Geodesic Rigidity]
    \label{cjc:totally_geodesic_rigidity}
    Suppose that $M$ is a compact Riemannian manifold with negative sectional curvature, of dimension at least $3$.
    Assume that for some $k\geq 2$ there exist infinitely many distinct immersed maximal totally geodesic $k$-dimensional submanifolds of $M$.

    Then $M$ is isometric to a negatively curved locally symmetric space, i.e. $\leftrightquot{\Lambda}{G}{K}$ where $G$ is a real rank $1$ Lie group, $K\subset G$ is a maximal compact, and $\Lambda$ is a lattice in $G$.
\end{conjecture}
Just like the result in \cite{FilipFisherLowe2024_Finiteness-of-totally-geodesic-hypersurfaces}, it should be possible to relax the assumptions to only require $M$ to have Anosov geodesic flow, and allow it to be noncompact but with appropriate bounds at infinity.
The conclusion of \cref{cjc:totally_geodesic_rigidity} is naturally strengthened by the work of Bader--Fisher--Miller--Stover \cite{BaderFisherMiller2021_Arithmeticity-superrigidity-and-totally-geodesic-submanifolds} to imply that the lattice $\Lambda$ defining $M$ is \emph{arithmetic}.

It is worth mentioning that a conclusion like in \cref{cjc:totally_geodesic_rigidity}, but assuming instead the smoothness of the stable and unstable foliations, was obtained by Benoist, Foulon, and Labourie
\cite{BenoistFoulonLabourie1992_Flots-dAnosov-a-distributions-stable-et-instable-differentiables}.
Zeghib \cite{Zeghib_Laminations-et-hypersurfaces-geodesiques-des-varietes1991} also considered closely related questions concerning totally geodesic manifolds.

\noindent \textbf{\Teichmuller Dynamics.}
In the setting of the action of $\SL_2(\bR)$ on moduli spaces of translation surfaces, measure and topological rigidity results were obtained by Eskin--Mirzakhani--Mohammadi \cite{EskinMirzakhani2018_Invariant-and-stationary-measures-for-the-rm-SL2R-action-on-moduli-space,EskinMirzakhaniMohammadi2015_Isolation-equidistribution-and-orbit-closures-for-the-rm-SL2R-action-on-moduli}.
The orbit closures that arise have a description entirely in terms of algebraic geometry \cite{Filip2016_Semisimplicity-and-rigidity-of-the-Kontsevich-Zorich-cocycle,Filip2016_Splitting-mixed-Hodge-structures-over-affine-invariant}, and furthermore the algebraic hull of the corresponding tangent cocycle satisfies very strong rigidity properties established in \cite{EskinFilipWright2018_The-algebraic-hull-of-the-Kontsevich-Zorich-cocycle}.
This, in turn allows us to establish a fairly complete version of \cref{pb:infinitely_many_implies_homogeneous}, in its ``finiteness'' form as explicated in \cref{rmk:on_homogeneity}.

\begin{remark}[On being special]
    \label{rmk:on_being_special}
    In the examples concerning totally geodesic submanifolds, and in \Teichmuller dynamics, the manifolds $Q$ that admit infinitely many maximal orbit closures carry, in fact, an additional structure.
    This structure then explains the presence of the infinitely many orbit closures.

    In the case of totally geodesic submanifolds, it is the arithmeticity of the lattice $\Lambda$, while in \Teichmuller dynamics the orbit closures are divided into ``typical'' and ``atypical'' (see \cite[\S5.1]{Filip2024_Translation-surfaces:-Dynamics-and-Hodge-theory}).
    The maximal atypical ones are always finite in number, and the typical ones are all explained by standard constructions.
\end{remark}

\noindent\textbf{Symplectic Manifolds.}
We end the discussion by addressing how to make sense of \cref{pb:infinitely_many_implies_homogeneous} in the context of symplectic manifolds and \cref{cjc:orbit_closures_on_symplectic_manifolds}.
The following class of homogeneous symplectic manifolds should be useful:
\begin{example}[Homogeneous symplectic manifolds]
    \label{eg:homogeneous_symplectic_manifolds}
    Suppose that $G$ is a simply connected Lie group acting transitively by symplectomorphisms on a simply connected symplectic manifold $(N,\omega)$.
    It is elementary to check that there exists an $\bR$-central extension $\widehat{G}\to G$, such that $(N,\omega)$ is equivariantly isomorphic to a universal cover of a coadjoint orbit of $G$ in $\widehat{\frakg}^*$.
    Note that $G$ acts on $\widehat{\frakg}^*$ since the coadjoint action of $\widehat{G}$ descends, as the center acts trivially.

    In the setting of \cref{cjc:orbit_closures_on_symplectic_manifolds}, one expects that if there are infinitely many maximal orbit closures $(N,\nu)$, then the ambient symplectic manifold $(Q,\omega)$ should be locally homogeneous, and hence locally isomorphic to a coadjoint orbit with the acting group $\Gamma_\mu$ contained in $G$.
\end{example}

%%                    end of subsec: Rigidity and Special Orbit Closures
%%=============================================================================

%%%%%%%%%%%%%%%%%%%%%%%%%%%%%%%%%%%%%%%%%%%%%%%%%%%%%%%%%%%%%%%%%%%%%%%%%%%%%%%%
%%%%% SubSection Dynamics on Character Varieties %%%%%

    \subsection{Dynamics on Character Varieties.}
        \label{ssec:Dynamics-on-Character-Varieties}

We discuss now a setting where we hope to have a classification of ergodic invariant measures and orbit closures.
Namely, we describe a rich class of subgroups of the mapping class group, namely those ``coming from algebraic geometry'', which act on the character variety of a compact Lie group.
For the action of these groups, one can prove results along the lines of Conjecture~\ref{cjc:orbit_of_a_point_on_a_k3_surface}.
The ergodic theory of the action of the full mapping class group on the character variety was initiated by Goldman \cite{Goldman1997_Ergodic-theory-on-moduli-spaces}.
This is joint work in preparation with Brown, Eskin, and Rodriguez Hertz.

To simplify notation, we will restrict to pure mapping class groups and character varieties of compact surfaces.
There is no additional difficulty in the proofs to extend to the case of relative character varieties and punctured surfaces.

\noindent\textbf{Character varieties of surface groups.}
Fix $g\geq 2$, let $\Sigma_g$ be a compact genus $g$ surface, and let $\pi_1(\Sigma_g)$ be its fundamental group (relative to a specified base point, omited from the notation).
Let also $U$ be a compact connected Lie group.
Introduce the \emph{representation variety}
\[
    \cR(g,U):=\Hom(\pi_1(\Sigma_g),U)
\]
which can be described as the preimage in $U^{2g}$ of the identity, under the map $(a_1,\dots,a_g,b_1,\dots,b_g)\mapsto [a_1,b_1]\dots[a_g,b_g]$.
Consider now the \emph{character variety}
\[
    \cX(g,U):={\cR(g,U)}\doublequot{U}
\]
where $U$ acts by conjugation.
Every compact Lie group admits a unique real algebraic structure, so we can consider the above quotients in an algebraic setting as well and equip the sets with algebraic structures.
In general, we get a space stratified by manifolds, and we can restrict our attention to the locus of representations which are Zariski-dense, which has a natural structure of smooth manifold.

\noindent\textbf{Mapping Class Groups and its Algebro-Geometric subgroups.}
Recall that $\Mod(\Sigma_g)$ denotes the mapping class group of $\Sigma_g$, which is the quotient $\Diff^+(\Sigma_g)/\Diff^0(\Sigma_g)$, i.e. the group of connected components of all orientation-preserving diffeomorphisms of $\Sigma_g$.
The mapping class group has two fundamental features: on the one hand it is the (orbifold) fundamental group of the moduli space of genus $g$ Riemann surfaces, and on the other it acts by outer automorphisms on $\pi_1(\Sigma_g)$, in fact the Dehn--Nielsen--Baer theorem \cite[Thm.~8.1]{FarbMargalit2012_A-primer-on-mapping-class-groups} identifies $\Mod(\Sigma_g)$ with an the index $2$ subgroup of $\operatorname{Out}(\pi_1(\Sigma_g))$ preserving orientation.
 % (viewed as an element of $H^2(\pi_1(\Sigma_g);\bZ)$).

Consider now a quasi-projective variety $B$, equipped with a projective family of Riemann surfaces $C\xrightarrow{\pi} B$.
If we pick a basepoint $b_0\in B$ and a diffeomorphism $\pi^{-1}(b_0)\isom \Sigma_g$, we obtain a monodromy representation $\rho\colon \pi_1(B,b)\to \Mod(\Sigma_g)$.

\begin{definition}[Algebro-Geometric subgroup of the Mapping Class Group]
    \label{def:algebro_geometric_subgroup_of_the_mapping_class_group}
    Call a subgroup $\Gamma\subset \Mod(\Sigma_g)$ \emph{algebro-geometric} if it arises as the image of a monodromy representation of a projective family of Riemann surfaces as above.
\end{definition}

\begin{remark}[On algebro-geometric subgroups]
    \label{rmk:on_algebro_geometric_subgroups}
    \leavevmode
    \begin{enumerate}
        \item The class of algebro-geometric subgroups is invariant under conjugation by elements of $\Mod(\Sigma_g)$, since part of the construction is the choice of an identification $\pi^{-1}(b_0)\isom \Sigma_g$.
        The full mapping class group $\Mod(\Sigma_g)$ is algebro-geometric.
        \item There is no loss of generality in assuming that $B$ is a Riemann surface of finite type, i.e. a compact Riemann surface with finitely many points removed.
        Indeed, starting from an arbitrary quasi-projective $B$, we can take a sufficiently general complete intersection curve in $B$ -- its fundamental group will surject onto that of $B$ and leads to the same algebro-geometric subgroup of $\Mod(\Sigma_g)$.
    \end{enumerate}
\end{remark}

To study the action of algebro-geometric subgroups of $\Mod(\Sigma_g)$ on $\cX(g,U)$ we can make use of Hodge theory.
Indeed, a Riemann surface structure on $\Sigma_g$ identifies $\cX(g,U)$ with a moduli space of semistable holomorphic bundles, and endows each tangent space at a point of the character variety with a weight $1$ Hodge structure.
This allows us to bring in the tools developed in \cite{Filip2016_Semisimplicity-and-rigidity-of-the-Kontsevich-Zorich-cocycle,Filip2017_Zero-Lyapunov-exponents-and-monodromy-of-the-Kontsevich-Zorich-cocycle}.
A key role is played by the subharmonicity of the Hodge norm, which is a replacement for uniform expansion in \cref{def:Uniform-Expansion-and-Uniform-Gaps}.

%%%%% End of SubSection Dynamics on Character Varieties %%%%%
%%%%%%%%%%%%%%%%%%%%%%%%%%%%%%%%%%%%%%%%%%%%%%%%%%%%%%%%%%%%%%%%%%%%%%%%%%%%%%%%

%%%%%%%%%%%%%%%%%%%%%%%%%%%%%%%%%%%%%%%%%%%%%%%%%%%%%%%%%%%%%%%%%%%%%%%%%%%%%%%%
%%%%%%%%%%%% End of Section Further Directions %%%%%%%%%%%%%%%%%%%%
%%%%%%%%%%%%%%%%%%%%%%%%%%%%%%%%%%%%%%%%%%%%%%%%%%%%%%%%%%%%%%%%%%%%%%%%%%%%%%%%

% acknowledgments
\section*{Acknowledgments.}
I am grateful to my collaborators on the works described here: Aaron Brown, Alex Eskin, David Fisher, Ben Lowe, and Federico Rodriguez Hertz.
I am grateful to Serge Cantat, Alex Eskin, and Federico Rodriguez Hertz for comments on a preliminary version of this text.
This material is based upon work supported by the National Science Foundation under Grants No. DMS-2005470 and DMS-2305394.
This research was partially conducted during the period the author served as a Clay Research Fellow.

% bibliography
\bibliographystyle{siamplain}
\bibliography{sfilip_rigidity_beyond_ICM}

\begin{thebibliography}{10}

\bibitem{AlvarezLeguilObata2022_Rigidity-of-U-Gibbs-measures-near-conservative-Anosov}
{\sc S.~Alvarez, M.~Leguil, D.~Obata, and B.~Santiago}, {\em Rigidity of
  {$\mathbf{\textit{U}}$-Gibbs} measures near conservative anosov
  diffeomorphisms on {$\mathbb{T}^3$}}, arXiv,  (2022).

\bibitem{BaderFisherMiller2021_Arithmeticity-superrigidity-and-totally-geodesic-submanifolds}
{\sc U.~Bader, D.~Fisher, N.~Miller, and M.~Stover}, {\em Arithmeticity,
  superrigidity, and totally geodesic submanifolds}, Ann. of Math. (2), 193
  (2021), pp.~837--861, \url{https://doi.org/10.4007/annals.2021.193.3.4},
  \url{https://doi.org/10.4007/annals.2021.193.3.4}.

\bibitem{BenoistFoulonLabourie1992_Flots-dAnosov-a-distributions-stable-et-instable-differentiables}
{\sc Y.~Benoist, P.~Foulon, and F.~c. Labourie}, {\em Flots d'{A}nosov \`a{}
  distributions stable et instable diff\'erentiables}, J. Amer. Math. Soc., 5
  (1992), pp.~33--74, \url{https://doi.org/10.2307/2152750},
  \url{https://doi.org/10.2307/2152750}.

\bibitem{BenoistQuint2011_Mesures-stationnaires-et-fermes-invariants-des-espaces-homogenes}
{\sc Y.~Benoist and J.-F. Quint}, {\em Mesures stationnaires et ferm\'{e}s
  invariants des espaces homog\`enes}, Ann. of Math. (2), 174 (2011),
  pp.~1111--1162, \url{https://doi.org/10.4007/annals.2011.174.2.8},
  \url{https://doi.org/10.4007/annals.2011.174.2.8}.

\bibitem{BenoistQuint2013_Stationary-measures-and-invariant-subsets-of-homogeneous-spaces-II}
{\sc Y.~Benoist and J.-F. Quint}, {\em Stationary measures and invariant
  subsets of homogeneous spaces ({II})}, J. Amer. Math. Soc., 26 (2013),
  pp.~659--734, \url{https://doi.org/10.1090/S0894-0347-2013-00760-2},
  \url{https://doi.org/10.1090/S0894-0347-2013-00760-2}.

\bibitem{BourgainFurmanLindenstrauss2011_Stationary-measures-and-equidistribution-for-orbits-of-nonabelian-semigroups-on-the-torus}
{\sc J.~Bourgain, A.~Furman, E.~Lindenstrauss, and S.~Mozes}, {\em Stationary
  measures and equidistribution for orbits of nonabelian semigroups on the
  torus}, J. Amer. Math. Soc., 24 (2011), pp.~231--280,
  \url{https://doi.org/10.1090/S0894-0347-2010-00674-1},
  \url{https://doi.org/10.1090/S0894-0347-2010-00674-1}.

\bibitem{Brown2023_Lattice-subgroups-acting-on-manifolds}
{\sc A.~Brown}, {\em Lattice subgroups acting on manifolds}, in
  I{CM}---{I}nternational {C}ongress of {M}athematicians. {V}ol. 5. {S}ections
  9--11, EMS Press, Berlin, 2023, pp.~3388--3411.

\bibitem{BrownEskinFilip2025_Measure-rigidity-for-generalized-u-Gibbs-states-and-stationary}
{\sc A.~Brown, A.~Eskin, S.~Filip, and F.~Rodriguez~Hertz}, {\em Measure
  rigidity for generalized u-{G}ibbs states and stationary measures via the
  factorization method}, arXiv:2502.14042,  (2025).

\bibitem{BrownRodriguez-Hertz2017_Measure-rigidity-for-random-dynamics-on-surfaces-and-related-skew}
{\sc A.~Brown and F.~Rodriguez~Hertz}, {\em Measure rigidity for random
  dynamics on surfaces and related skew products}, J. Amer. Math. Soc., 30
  (2017), pp.~1055--1132, \url{https://doi.org/10.1090/jams/877},
  \url{https://doi.org/10.1090/jams/877}.

\bibitem{CantatDujardin2023_Dynamics-of-automorphism-groups-of-projective-surfaces:-classification-examples}
{\sc S.~Cantat and R.~Dujardin}, {\em Dynamics of automorphism groups of
  projective surfaces: classification, examples and outlook}, 2023,
  \url{https://doi.org/10.48550/ARXIV.2310.01303},
  \url{https://arxiv.org/abs/2310.01303}.

\bibitem{CantatDujardin2023_Random-dynamics-on-real-and-complex-projective-surfaces}
{\sc S.~Cantat and R.~Dujardin}, {\em Random dynamics on real and complex
  projective surfaces}, J. Reine Angew. Math., 802 (2023), pp.~1--76,
  \url{https://doi.org/10.1515/crelle-2023-0038}.

\bibitem{CantatDujardin2024_Finite-orbits-for-large-groups-of-automorphisms-of-projective}
{\sc S.~Cantat and R.~Dujardin}, {\em Finite orbits for large groups of
  automorphisms of projective surfaces}, Compos. Math., 160 (2024),
  pp.~120--175, \url{https://doi.org/10.1112/s0010437x23007613},
  \url{https://doi.org/10.1112/s0010437x23007613}.

\bibitem{CantatDujardin2025_Hyperbolicity-for-large-automorphism-groups-of-projective-surfaces}
{\sc S.~Cantat and R.~Dujardin}, {\em Hyperbolicity for large automorphism
  groups of projective surfaces}, J. \'Ec. polytech. Math., 12 (2025),
  pp.~421--480, \url{https://doi.org/10.5802/jep.294},
  \url{https://doi.org/10.5802/jep.294}.

\bibitem{CantatDupont2020_Automorphisms-of-surfaces:-Kummer-rigidity-and-measure-of-maximal-entropy}
{\sc S.~Cantat and C.~Dupont}, {\em Automorphisms of surfaces: {K}ummer
  rigidity and measure of maximal entropy}, J. Eur. Math. Soc. (JEMS), 22
  (2020), pp.~1289--1351, \url{https://doi.org/10.4171/JEMS/946},
  \url{https://doi.org/10.4171/JEMS/946}.

\bibitem{Chung2020_Stationary-measures-and-orbit-closures-of-uniformly-expanding-random}
{\sc P.~N. Chung}, {\em Stationary measures and orbit closures of uniformly
  expanding random dynamical systems on surfaces},  (2020),
  \url{https://doi.org/https://doi.org/10.48550/arXiv.2006.03166},
  \url{https://doi.org/10.48550/arXiv.2006.03166}.

\bibitem{Crauel1993_Non-Markovian-invariant-measures-are-hyperbolic}
{\sc H.~Crauel}, {\em Non-{M}arkovian invariant measures are hyperbolic},
  Stochastic Process. Appl., 45 (1993), pp.~13--28,
  \url{https://doi.org/10.1016/0304-4149(93)90057-B},
  \url{https://doi.org/10.1016/0304-4149(93)90057-B}.

\bibitem{MR1101994}
{\sc S.~G. Dani and G.~A. Margulis}, {\em Asymptotic behaviour of trajectories
  of unipotent flows on homogeneous spaces}, Proc. Indian Acad. Sci. Math.
  Sci., 101 (1991), pp.~1--17, \url{https://doi.org/10.1007/BF02872005},
  \url{https://doi.org/10.1007/BF02872005}.

\bibitem{DaniMargulis1993_Limit-distributions-of-orbits-of-unipotent-flows-and-values}
{\sc S.~G. Dani and G.~A. Margulis}, {\em Limit distributions of orbits of
  unipotent flows and values of quadratic forms}, in I. {M}. {G}elfand
  {S}eminar, vol.~16, Part 1 of Adv. Soviet Math., Amer. Math. Soc.,
  Providence, RI, 1993, pp.~91--137.

\bibitem{DeWittDolgopyat2025_Conservative-Coexpanding-on-Average-Diffeomorphisms}
{\sc J.~DeWitt and D.~Dolgopyat}, {\em Conservative coexpanding on average
  diffeomorphisms}, arXiv,  (2025),
  \url{https://doi.org/10.48550/arXiv.2503.06855}.

\bibitem{DolgopyatKrikorian2007_On-simultaneous-linearization-of-diffeomorphisms-of-the-sphere}
{\sc D.~Dolgopyat and R.~Krikorian}, {\em On simultaneous linearization of
  diffeomorphisms of the sphere}, Duke Math. J., 136 (2007), pp.~475--505.

\bibitem{Elliott-Smith2023_Uniformly-expanding-random-walks-on-manifolds}
{\sc R.~Elliott~Smith}, {\em Uniformly expanding random walks on manifolds},
  Nonlinearity, 36 (2023), pp.~5955--5972,
  \url{https://doi.org/10.1088/1361-6544/acfa5a},
  \url{https://doi.org/10.1088/1361-6544/acfa5a}.

\bibitem{EskinFilipWright2018_The-algebraic-hull-of-the-Kontsevich-Zorich-cocycle}
{\sc A.~Eskin, S.~Filip, and A.~Wright}, {\em The algebraic hull of the
  {K}ontsevich-{Z}orich cocycle}, Ann. of Math. (2), 188 (2018), pp.~281--313,
  \url{https://doi.org/10.4007/annals.2018.188.1.5},
  \url{https://doi.org/10.4007/annals.2018.188.1.5}.

\bibitem{EskinLindenstrauss2018_Random-walks-on-locally-homogeneous-spaces}
{\sc A.~Eskin and E.~Lindenstrauss}, {\em Random walks on locally homogeneous
  spaces}.
\newblock preprint, 2018,
  \url{https://math.uchicago.edu/~eskin/RandomWalks/paper.pdf}.

\bibitem{EskinMirzakhani2018_Invariant-and-stationary-measures-for-the-rm-SL2R-action-on-moduli-space}
{\sc A.~Eskin and M.~Mirzakhani}, {\em Invariant and stationary measures for
  the {${\rm SL}(2,\mathbb{R})$} action on moduli space}, Publ. Math. Inst.
  Hautes \'{E}tudes Sci., 127 (2018), pp.~95--324,
  \url{https://doi.org/10.1007/s10240-018-0099-2},
  \url{https://doi.org/10.1007/s10240-018-0099-2}.

\bibitem{EskinMirzakhaniMohammadi2015_Isolation-equidistribution-and-orbit-closures-for-the-rm-SL2R-action-on-moduli}
{\sc A.~Eskin, M.~Mirzakhani, and A.~Mohammadi}, {\em Isolation,
  equidistribution, and orbit closures for the {${\rm SL}(2,\mathbb{R})$}
  action on moduli space}, Ann. of Math. (2), 182 (2015), pp.~673--721,
  \url{https://doi.org/10.4007/annals.2015.182.2.7},
  \url{https://doi.org/10.4007/annals.2015.182.2.7}.

\bibitem{EskinPotrieZhang2023_Geometric-properties-of-partially-hyperbolic-measures-and-applications}
{\sc A.~Eskin, R.~Potrie, and Z.~Zhang}, {\em Geometric properties of partially
  hyperbolic measures and applications to measure rigidity}, 2023,
  \url{https://arxiv.org/abs/2302.12981}.

\bibitem{FarbMargalit2012_A-primer-on-mapping-class-groups}
{\sc B.~Farb and D.~Margalit}, {\em A primer on mapping class groups}, vol.~49
  of Princeton Mathematical Series, Princeton University Press, Princeton, NJ,
  2012.

\bibitem{Filip2016_Semisimplicity-and-rigidity-of-the-Kontsevich-Zorich-cocycle}
{\sc S.~Filip}, {\em Semisimplicity and rigidity of the {K}ontsevich-{Z}orich
  cocycle}, Invent. Math., 205 (2016), pp.~617--670,
  \url{https://doi.org/10.1007/s00222-015-0643-3},
  \url{https://doi.org/10.1007/s00222-015-0643-3}.

\bibitem{Filip2016_Splitting-mixed-Hodge-structures-over-affine-invariant}
{\sc S.~Filip}, {\em Splitting mixed {H}odge structures over affine invariant
  manifolds}, Ann. of Math. (2), 183 (2016), pp.~681--713,
  \url{https://doi.org/10.4007/annals.2016.183.2.5},
  \url{https://doi.org/10.4007/annals.2016.183.2.5}.

\bibitem{Filip2017_Zero-Lyapunov-exponents-and-monodromy-of-the-Kontsevich-Zorich-cocycle}
{\sc S.~Filip}, {\em Zero {L}yapunov exponents and monodromy of the
  {K}ontsevich-{Z}orich cocycle}, Duke Math. J., 166 (2017), pp.~657--706,
  \url{https://doi.org/10.1215/00127094-3715806},
  \url{https://doi.org/10.1215/00127094-3715806}.

\bibitem{Filip2022_An-introduction-to-K3-surfaces-and-their-dynamics}
{\sc S.~Filip}, {\em An introduction to {K}3 surfaces and their dynamics}, in
  Teichm\"uller theory and dynamics, vol.~58 of Panor. Synth\`eses, Soc. Math.
  France, Paris, 2022, pp.~1--46.

\bibitem{Filip2024_Translation-surfaces:-Dynamics-and-Hodge-theory}
{\sc S.~Filip}, {\em Translation surfaces: Dynamics and hodge theory}, EMS
  Surveys in Mathematical Sciences, 11 (2024), pp.~63--151,
  \url{https://doi.org/10.4171/emss/78},
  \url{http://dx.doi.org/10.4171/EMSS/78}.

\bibitem{FilipFisherLowe2024_Finiteness-of-totally-geodesic-hypersurfaces}
{\sc S.~Filip, D.~Fisher, and B.~Lowe}, {\em Finiteness of totally geodesic
  hypersurfaces}, arXiv:2408.03430,  (2024),
  \url{https://doi.org/10.48550/arXiv.2408.03430}.

\bibitem{FilipTosatti2021_Kummer-rigidity-for-K3-surface-automorphisms-via-Ricci-flat-metrics}
{\sc S.~Filip and V.~Tosatti}, {\em Kummer rigidity for {K3} surface
  automorphisms via {R}icci-flat metrics}, Amer. J. Math., 143 (2021),
  pp.~1431--1462, \url{https://doi.org/10.1353/ajm.2021.0036},
  \url{https://doi.org/10.1353/ajm.2021.0036}.

\bibitem{FilipTosatti2023_Canonical-currents-and-heights-for-K3-surfaces}
{\sc S.~Filip and V.~Tosatti}, {\em Canonical currents and heights for {K}3
  surfaces}, Camb. J. Math., 11 (2023), pp.~699--794,
  \url{https://doi.org/10.4310/cjm.2023.v11.n3.a2},
  \url{https://doi.org/10.4310/cjm.2023.v11.n3.a2}.

\bibitem{Fisher2023_Rigidity-lattices-and-invariant-measures-beyond-homogeneous}
{\sc D.~Fisher}, {\em Rigidity, lattices, and invariant measures beyond
  homogeneous dynamics}, in I{CM}---{I}nternational {C}ongress of
  {M}athematicians. {V}ol. 5. {S}ections 9--11, EMS Press, Berlin, 2023,
  pp.~3484--3507.

\bibitem{Furman2002_Random-walks-on-groups-and-random-transformations}
{\sc A.~Furman}, {\em Random walks on groups and random transformations}, in
  Handbook of dynamical systems, {V}ol. 1{A}, North-Holland, Amsterdam, 2002,
  pp.~931--1014, \url{https://doi.org/10.1016/S1874-575X(02)80014-5},
  \url{https://doi.org/10.1016/S1874-575X(02)80014-5}.

\bibitem{Furstenberg1963_Noncommuting-random-products}
{\sc H.~Furstenberg}, {\em Noncommuting random products}, Trans. Amer. Math.
  Soc., 108 (1963), pp.~377--428, \url{https://doi.org/10.2307/1993589},
  \url{https://doi.org/10.2307/1993589}.

\bibitem{Goldman1997_Ergodic-theory-on-moduli-spaces}
{\sc W.~M. Goldman}, {\em Ergodic theory on moduli spaces}, Ann. of Math. (2),
  146 (1997), pp.~475--507, \url{https://doi.org/10.2307/2952454},
  \url{https://doi.org/10.2307/2952454}.

\bibitem{Gromov1988_Rigid-transformations-groups}
{\sc M.~Gromov}, {\em Rigid transformations groups}, in G\'eom\'etrie
  diff\'erentielle ({P}aris, 1986), vol.~33 of Travaux en Cours, Hermann,
  Paris, 1988, pp.~65--139.

\bibitem{GuysinskyKatok1998_Normal-forms-and-invariant-geometric-structures-for-dynamical-systems}
{\sc M.~Guysinsky and A.~Katok}, {\em Normal forms and invariant geometric
  structures for dynamical systems with invariant contracting foliations},
  Math. Res. Lett., 5 (1998), pp.~149--163,
  \url{https://doi.org/10.4310/MRL.1998.v5.n2.a2},
  \url{https://doi.org/10.4310/MRL.1998.v5.n2.a2}.

\bibitem{KalininSadovskaya2017_Normal-forms-for-non-uniform-contractions}
{\sc B.~Kalinin and V.~Sadovskaya}, {\em Normal forms for non-uniform
  contractions}, J. Mod. Dyn., 11 (2017), pp.~341--368,
  \url{https://doi.org/10.3934/jmd.2017014},
  \url{https://doi.org/10.3934/jmd.2017014}.

\bibitem{Katz2023_Measure-rigidity-of-Anosov-flows-via-the-factorization-method}
{\sc A.~Katz}, {\em Measure rigidity of {A}nosov flows via the factorization
  method}, Geom. Funct. Anal., 33 (2023), pp.~468--540,
  \url{https://doi.org/10.1007/s00039-023-00629-8},
  \url{https://doi.org/10.1007/s00039-023-00629-8}.

\bibitem{Ledrappier1986_Positivity-of-the-exponent-for-stationary-sequences-of-matrices}
{\sc F.~Ledrappier}, {\em Positivity of the exponent for stationary sequences
  of matrices}, in Lyapunov exponents ({B}remen, 1984), vol.~1186 of Lecture
  Notes in Math., Springer, Berlin, 1986, pp.~56--73,
  \url{https://doi.org/10.1007/BFb0076833},
  \url{https://doi.org/10.1007/BFb0076833}.

\bibitem{LedrappierYoung1985_The-metric-entropy-of-diffeomorphisms.-II.-Relations-between-entropy-exponents}
{\sc F.~Ledrappier and L.-S. Young}, {\em The metric entropy of
  diffeomorphisms. {II}. {R}elations between entropy, exponents and dimension},
  Ann. of Math. (2), 122 (1985), pp.~540--574,
  \url{https://doi.org/10.2307/1971329}, \url{https://doi.org/10.2307/1971329}.

\bibitem{Lindenstrauss2010_Equidistribution-in-homogeneous-spaces-and-number-theory}
{\sc E.~Lindenstrauss}, {\em Equidistribution in homogeneous spaces and number
  theory}, in Proceedings of the {I}nternational {C}ongress of
  {M}athematicians. {V}olume {I}, Hindustan Book Agency, New Delhi, 2010,
  pp.~531--557.

\bibitem{Liu2016_Lyapunov-Exponents-Approximation-Symplectic-Cocycle-Deformation}
{\sc X.-C. Liu}, {\em Lyapunov Exponents Approximation, Symplectic Cocycle
  Deformation and a Large Deviation Theorem}, PhD thesis, IMPA, 2016,
  \url{https://doi.org/https://impa.br/wp-content/uploads/2017/11/tese_dout_Xiaochuan_Liu.pd}.

\bibitem{MohammadiMargulis2022_Arithmeticity-of-hyperbolic-3-manifolds-containing-infinitely-many-totally}
{\sc G.~Margulis and A.~Mohammadi}, {\em Arithmeticity of hyperbolic
  3-manifolds containing infinitely many totally geodesic surfaces}, Ergodic
  Theory Dynam. Systems, 42 (2022), pp.~1188--1219,
  \url{https://doi.org/10.1017/etds.2021.21},
  \url{https://doi.org/10.1017/etds.2021.21}.

\bibitem{Margulis1987_Formes-quadratriques-indefinies-et-flots-unipotents-sur-les-espaces-homogenes}
{\sc G.~A. Margulis}, {\em Formes quadratriques ind\'{e}finies et flots
  unipotents sur les espaces homog\`enes}, C. R. Acad. Sci. Paris S\'{e}r. I
  Math., 304 (1987), pp.~249--253.

\bibitem{Margulis_Discrete-subgroups-of-semisimple-Lie-groups1991}
{\sc G.~A. Margulis}, {\em Discrete subgroups of semisimple {L}ie groups},
  vol.~17 of Ergebnisse der Mathematik und ihrer Grenzgebiete (3) [Results in
  Mathematics and Related Areas (3)], Springer-Verlag, Berlin, 1991,
  \url{https://doi.org/10.1007/978-3-642-51445-6},
  \url{https://doi.org/10.1007/978-3-642-51445-6}.

\bibitem{MargulisTomanov1994_Invariant-measures-for-actions-of-unipotent-groups-over-local}
{\sc G.~A. Margulis and G.~M. Tomanov}, {\em Invariant measures for actions of
  unipotent groups over local fields on homogeneous spaces}, Invent. Math., 116
  (1994), pp.~347--392, \url{https://doi.org/10.1007/BF01231565},
  \url{https://doi.org/10.1007/BF01231565}.

\bibitem{MR2958928}
{\sc J.~Mather}, {\em Notes on topological stability}, Bull. Amer. Math. Soc.
  (N.S.), 49 (2012), pp.~475--506,
  \url{https://doi.org/10.1090/S0273-0979-2012-01383-6},
  \url{https://doi.org/10.1090/S0273-0979-2012-01383-6}.

\bibitem{Melnick_Non-stationary-smooth-geometric-structures-for-contracting-measurable}
{\sc K.~Melnick}, {\em Non-stationary smooth geometric structures for
  contracting measurable cocycles}, Ergodic Theory Dynam. Systems, 39 (2019),
  pp.~392--424, \url{https://doi.org/10.1017/etds.2017.38},
  \url{https://doi.org/10.1017/etds.2017.38}.

\bibitem{PesinSinai1982_Gibbs-measures-for-partially-hyperbolic-attractors}
{\sc Y.~B. Pesin and Y.~G. Sinai}, {\em Gibbs measures for partially hyperbolic
  attractors}, Ergodic Theory Dynam. Systems, 2 (1982), pp.~417--438,
  \url{https://doi.org/10.1017/S014338570000170X},
  \url{https://doi.org/10.1017/S014338570000170X}.

\bibitem{Ratner1990_On-measure-rigidity-of-unipotent-subgroups-of-semisimple-groups}
{\sc M.~Ratner}, {\em On measure rigidity of unipotent subgroups of semisimple
  groups}, Acta Math., 165 (1990), pp.~229--309,
  \url{https://doi.org/10.1007/BF02391906},
  \url{https://doi.org/10.1007/BF02391906}.

\bibitem{Ratner1991_On-Raghunathans-measure-conjecture}
{\sc M.~Ratner}, {\em On {R}aghunathan's measure conjecture}, Ann. of Math.
  (2), 134 (1991), pp.~545--607, \url{https://doi.org/10.2307/2944357},
  \url{https://doi.org/10.2307/2944357}.

\bibitem{Ratner1991_Raghunathans-topological-conjecture-and-distributions-of-unipotent-flows}
{\sc M.~Ratner}, {\em Raghunathan's topological conjecture and distributions of
  unipotent flows}, Duke Math. J., 63 (1991), pp.~235--280,
  \url{https://doi.org/10.1215/S0012-7094-91-06311-8},
  \url{https://doi.org/10.1215/S0012-7094-91-06311-8}.

\bibitem{Roda2024_Classifying-hyperbolic-ergodic-stationary-measures-on-compact}
{\sc M.~Roda}, {\em Classifying hyperbolic ergodic stationary measures on
  compact complex surfaces with large automorphism groups}, arXiv,  (2024),
  \url{https://doi.org/10.48550/arXiv.2410.18350}.

\bibitem{Sharpe1997_Differential-geometry}
{\sc R.~W. Sharpe}, {\em Differential geometry}, vol.~166 of Graduate Texts in
  Mathematics, Springer-Verlag, New York, 1997.

\bibitem{Sternberg1957_Local-contractions-and-a-theorem-of-Poincare}
{\sc S.~Sternberg}, {\em Local contractions and a theorem of {P}oincar\'{e}},
  Amer. J. Math., 79 (1957), pp.~809--824,
  \url{https://doi.org/10.2307/2372437}, \url{https://doi.org/10.2307/2372437}.

\bibitem{Zeghib_Laminations-et-hypersurfaces-geodesiques-des-varietes1991}
{\sc A.~Zeghib}, {\em Laminations et hypersurfaces g\'{e}od\'{e}siques des
  vari\'{e}t\'{e}s hyperboliques}, Ann. Sci. \'{E}cole Norm. Sup. (4), 24
  (1991), pp.~171--188,
  \url{http://www.numdam.org/item?id=ASENS_1991_4_24_2_171_0}.

\bibitem{Zimmer1987_Actions-of-semisimple-groups-and-discrete-subgroups}
{\sc R.~J. Zimmer}, {\em Actions of semisimple groups and discrete subgroups},
  in Proceedings of the {I}nternational {C}ongress of {M}athematicians, {V}ol.
  1, 2 ({B}erkeley, {C}alif., 1986), Amer. Math. Soc., Providence, RI, 1987,
  pp.~1247--1258.

\end{thebibliography}

\crefname{subsection}{section}{Sections}

\end{document}